%% file: root.tex
\documentclass[twoside]{article} 

\usepackage{jmlr2e} 
\usepackage{graphicx} 
\usepackage{subfigure} 

\usepackage{natbib}

\usepackage{algorithm}
\usepackage{algorithmic}

\usepackage{hyperref}

\usepackage{amsmath}
\usepackage{pgfplots}
\usepackage{color}
\input{my_sections.tex}

\usepackage{flushend}
\usepackage{multirow}
\usepackage{rotating}
\usepackage{appendix}
\usepackage{tikz,pgfplots}

\usepackage{standalone} 

\newif\ifcompiletikz
\compiletikztrue

\input{mysymbol.sty}

\newtheorem{assumption}{\hspace{0pt}\bf Assumption}

\begin{document}


\jmlrheading{1}{2000}{1-48}{4/00}{10/00}{Aryan Mokhtari and Alejandro Ribeiro}


\ShortHeadings{DSA: Decentralized Double Stochastic Averaging Gradient Algorithm}{Mokhtari and Ribeiro}
\firstpageno{1}



\title{DSA: Decentralized Double Stochastic Averaging Gradient Algorithm}

\author{\name Aryan Mokhtari  \email aryanm@seas.upenn.edu \\
       \addr Department of Electrical and Systems Engineering\\
       University of Pennsylvania\\
       Philadelphia, PA 19104, USA
       \AND
       \name Alejandro Ribeiro \email aribeiro@seas.upenn.edu  \\
       \addr Department of Electrical and Systems Engineering\\
       University of Pennsylvania\\
       Philadelphia, PA 19104, USA}

\editor{}

\maketitle


\thispagestyle{empty}

\begin{abstract}
This paper considers convex optimization problems where nodes of a network have access to summands of a global objective. Each of these local objectives is further assumed to be an average of a finite set of functions. The motivation for this setup is to solve large scale machine learning problems where elements of the training set are distributed to multiple computational elements. The decentralized double stochastic averaging gradient (DSA) algorithm is proposed as a solution alternative that relies on: (i) The use of local stochastic averaging gradients. (ii) Determination of descent steps as differences of consecutive stochastic averaging gradients. Strong convexity of local functions and Lipschitz continuity of local gradients is shown to guarantee linear convergence of the sequence generated by DSA in expectation. Local iterates are further shown to approach the optimal argument for almost all realizations. The expected linear convergence of DSA is in contrast to the sublinear rate characteristic of existing methods for decentralized stochastic optimization. Numerical experiments on a logistic regression problem illustrate reductions in convergence time and number of feature vectors processed until convergence relative to these other alternatives.

\end{abstract}


\input{Introduction.tex}

\input{Problem-1.tex}

\input{saddle_point_interpretation.tex}

\input{Convergence2.tex}

\input{Simulations.tex}

\input{Conclusions.tex}

\acks{We acknowledge the support of the National Science Foundation (NSF CAREER
CCF-0952867) and the Office of Naval Research (ONR N00014-12-1-0997).}

\input{Appendix.tex}
\bibliography{bmc_article}
\bibliographystyle{icml2014}
  
\end{document}

%% file: my_sections.tex
\usepackage{needspace}



%% file: Introduction.tex

\section{Introduction}\label{sec_Introduction}

We consider machine learning problems with large training sets that are distributed into a network of computing agents so that each of the nodes maintains a moderate number of samples. This leads to decentralized consensus optimization problems where summands of the global objective function are available at different nodes of the network. In this class of problems agents (nodes) try to optimize the global cost function by operating on their local functions and communicating with their neighbors only. Specifically, consider a variable $\bbx\in \reals^p$ and a connected network of size $N$ where each node $n$ has access to a local objective function $f_{n}: \reals^p \to \reals$. The local objective function $f_n(\bbx)$ is defined as the average of $q_n$ local instantaneous functions $f_{n,i}(\bbx)$ that can be individually evaluated at node $n$. Agents cooperate to solve the global optimization 
\begin{equation}\label{opt_problem}
   {\tbx^*}  := \argmin_\bbx \sum_{n=1}^N  f_{n}(\bbx)
   	   	 = \argmin_\bbx  \sum_{n=1}^N \frac{1}{q_n} \sum_{i=1}^{q_n}f_{n,i}(\bbx).
\end{equation} 
The formulation in \eqref{opt_problem} models a training set with a total of $\sum_{n=1}^{N}q_n$ training samples that are distributed among the $N$ agents for parallel processing conducive to the determination of the optimal classifier $\tbx^*$ (\cite{bekkerman2011scaling,Tsianos2012-allerton-consensus,Cevher2014}). Although we make no formal assumption, in cases of practical importance the total number of training samples $\sum_{n=1}^{N}q_n$ is very large, but the numbers of elements $q_n$ available at a specific node are moderate.

Our interest here is in solving \eqref{opt_problem} with a method that is decentralized -- nodes operate on their local functions and communicate with neighbors only  --, stochastic -- nodes determine a descent direction by evaluating only one out of the $q_n$ functions $f_{n,i}$ at each iteration --, and has a linear convergence rate in expectation -- the expected distance to the optimum is scaled by a subunit factor at each iteration. 

Decentralized optimization is relatively mature and various methods are known with complementary advantages. These methods include decentralized gradient descent (DGD) (\cite{Nedic2009,Jakovetic2014-1,YuanQing}), network Newton (\cite{mokhtari2015network,mokhtari2015network2}), decentralized dual averaging (\cite{Duchi2012,cTsianosEtal12}), the exact first order algorithm (EXTRA) (\cite{Shi2014}), as well as the alternating direction method of multipliers (ADMM) (\cite{BoydEtalADMM11,Shi2014-ADMM,iutzeler2013explicit}) and its linearized variants (\cite{cQingRibeiroADMM14,ling2014dlm,DQMglobalSIP}). The ADMM, its variants, and EXTRA converge linearly to the optimal argument but DGD, network Newton, and decentralized dual averaging have sublinear convergence rates. Of particular importance to this paper, is the fact that DGD has (inexact) linear converge to a neighborhood of the optimal argument when it uses constant stepsizes. It can achieve exact convergence by using diminishing stepsizes, but the convergence rate degrades to sublinear. This lack of linear convergence is solved by EXTRA through the use of iterations that rely on information of two consecutive steps (\cite{Shi2014}). 

All of the algorithms mentioned above require the computationally costly evaluation of the local gradients $\nabla f_n(\bbx)=(1/q_n)\sum_{i=1}^{q_n} \nabla f_{n,i}(\bbx)$. This cost can be avoided by stochastic decentralized algorithms that reduce computational cost of iterations by  substituting all local gradients with their stochastic approximations. 
This reduces the computational cost per iteration but results in sublinear convergence rates of order $O(1/t)$ even if the corresponding deterministic algorithm exhibits linear convergence. This is a drawback that also exists in centralized stochastic optimization where linear convergence rates in expectation are established by decreasing the variance of the stochastic gradient approximation (\cite{roux2012stochastic,schmidt2013minimizing,shalev2013stochastic,johnson2013accelerating,konevcny2013semi,defazio2014saga}). In this paper we build on the ideas of the stochastic averaging gradient (SAG) algorithm (\cite{schmidt2013minimizing}) and its unbiased version SAGA (\cite{defazio2014saga}). Both of these algorithms use the idea of stochastic incremental averaging gradients. At each iteration only one of the stochastic gradients is updated and the average of all of the most recent stochastic gradients is used for estimating gradient.  

The contribution of this paper is to develop the decentralized double stochastic averaging gradient (DSA) method, a novel decentralized stochastic algorithm for solving  \eqref{opt_problem}. The method exploits a new interpretation of EXTRA as a saddle point method and uses stochastic averaging gradients in lieu of gradients. DSA is decentralized because it is implementable in a network setting where nodes can communicate only with their neighbors. It is double because iterations utilize the information of two consecutive iterates. It is stochastic because the gradient of only one randomly selected function is evaluated at each iteration and it is an averaging method because it uses an average of stochastic gradients to approximate the local gradients. DSA is proven to converge linearly to the optimal argument $\tbx^*$ in expectation. This is in contrast to all other decentralized stochastic methods to solve \eqref{opt_problem} that converge at sublinear rates. 

We begin the paper with a discussion of DGD, EXTRA and stochastic averaging gradient. With these definitions in place we define the DSA algorithm by replacing the gradients used in EXTRA by stochastic averaging gradients (Section \ref{sec:DSA}). We follow with a digression on the limit points of DGD and EXTRA iterations to explain the reason why DGD does not achieve exact convergence but EXTRA is expected to do so (Section \ref{subset:lmit_point}). A reinterpretation of EXTRA as a saddle point method that solves for the critical points of the augmented Lagrangian of a constrained optimization problem equivalent to \eqref{opt_problem} is then introduced. It follows from this reinterpretation that DSA is a stochastic saddle point method (Section \ref{subset:saddle}). 
The fact that DSA is a stochastic saddle point method is the critical enabler of the subsequent convergence analysis (Section \ref{sec_convergence}). In particular, it is possible to guarantee that strong convexity and gradient Lipschitz continuity of the local instantaneous functions $f_{n,i}$ imply that a Lyapunov function associated with the sequence of iterates generated by DSA converges linearly to its optimal value in expectation (Theorem \ref{lin_convg_thm}). Linear convergence in expectation of the local iterates to the optimal argument $\tbx^*$ of \eqref{opt_problem} follows as a trivial consequence (Corollary \ref{corollary_lin_convg}). We complement this result by showing convergence of all the local variables to the optimal argument $\tbx^*$ with probability 1 (Theorem \ref{a_s_convg}).

The advantages of DSA relative to a group of stochastic and deterministic alternatives in solving a logistic regression problem with a synthetic dataset are then studied in numerical experiments (Section \ref{sec:simulations}). These results demonstrate that DSA is the only decentralized stochastic algorithm that reaches the optimal solution with a linear convergence rate. We further show that DSA outperforms deterministic algorithms when the metric is the number of times that elements of the training set are evaluated. The behavior of DSA for different network topologies is also evaluated. We close the paper with pertinent remarks (Section \ref{sec_conclusions}).

\medskip\noindent{\bf Notation\quad} Lowercase boldface $\bbv$ denotes a vector and uppercase boldface $\bbA$ a matrix. For column vectors $\bbx_1,\dots,\bbx_N$ we use the notation $\bbx=[\bbx_1;\ldots;\bbx_N]$ to represent the stack column vector $\bbx$. We use $\|\bbv\|$ to denote the Euclidean norm of vector $\bbv$ and $\|\bbA\|$ to denote the Euclidean norm of matrix $\bbA$. For a vector $\bbv$ and a positive definite matrix $\bbA$, the $\bbA$-weighted norm is defined as $\|\bbv\|_{\bbA}:=\sqrt{\bbv^T\bbA\bbv}$. The null space of matrix $\bbA$ is denoted by $\text{null}(\bbA)$ and the span of a vector by $\text{span}(\bbx)$. The operator $\mbE_{\bbx}[\cdot]$ stands for expectation over random variable $\bbx$ and  $\mbE[\cdot]$ for expectation with respect to the distribution of a stochastic process.

%% file: Problem-1.tex

%
\section{Decentralized Double stochastic averaging gradient }\label{sec:DSA}

Consider a connected network that contains $N$ nodes such that each node $n$ can only communicate with peers in its neighborhood $\ccalN_{n}$. Define $\bbx_n\in \reals^p$ as a local copy of the variable $\bbx$ that is kept at node $n$. In decentralized optimization, agents try to minimize their local functions $f_{n}(\bbx_n)$ while ensuring that their local variables $\bbx_n$ coincide with the variables $\bbx_{m}$ of all neighbors $m \in \ccalN_{n}$ -- which, given that the network is connected, ensures that the variables $\bbx_n$ of all nodes are the same and renders the problem equivalent to \eqref{opt_problem}. DGD is a well known method for decentralized optimization that relies on the introduction of nonnegative weights $w_{ij}\geq0$ that are not null if and only if  $m=n$ or if $m\in \mathcal{N}_n$. Letting $t\in\naturals$ be a discrete time index and $\alpha$ a given stepsize, DGD is defined by the recursion
\begin{equation}\label{gd_iteration}
   \bbx_{n}^{t+1} 
      = \sum_{m=1}^{N} w_{nm}\bbx_{m}^{t}-\alpha\nabla f_n({\bbx_{n}^{t}}), 
      \qquad n=1,\ldots,N.
\end{equation}
Since $w_{nm}=0$ when $m\neq n$ and $m\notin \mathcal{N}_n$, it follows from \eqref{gd_iteration} that node $n$ updates $\bbx_n$ by performing an average over the variables $\bbx_{m}^t$ of its neighbors $m\in \mathcal{N}_n$ and its own $\bbx_{n}^t$, followed by descent through the negative local gradient $-\nabla f_n(\bbx_{n}^t)$. If a constant stepsize is used, DGD iterates $\bbx_{n}^t$ approach a neighborhood of the optimal argument $\tbx^*$ of \eqref{opt_problem} but don't converge exactly. To achieve exact convergence diminishing stepsizes are used but the resulting convergence rate is sublinear (\cite{Nedic2009}). 

EXTRA is a method that resolves either of these issues by mixing two consecutive DGD iterations with different weight matrices and opposite signs. To be precise, introduce a second set of weights $\tdw_{nm}$ with the same properties as the weights $w_{nm}$ and define EXTRA through the recursion
\begin{align}\label{local_extra}
    \bbx_n^{t+1} = \bbx_n^{t}
                    + \sum_{m=1}^N w_{nm} \bbx_{m}^t 
                     - \sum_{m=1}^N \tdw_{nm}\bbx_{m}^{t-1}  - \alpha\left [ \nabla f_n(\bbx_n^{t}) 
                                    -\nabla f_n(\bbx_n^{t-1}) \right] , 
                     \quad  n=1,\ldots,N.
\end{align}
Observe that \eqref{local_extra} is well defined for $t>0$. For $t=0$ we utilize the regular DGD iteration in \eqref{gd_iteration}. In the nomenclature of this paper we say that EXTRA performs a decentralized double gradient descent step because it operates in a decentralized manner while utilizing a difference of two gradients as descent direction. Minor modification as it is, the use of this gradient difference in lieu of simple gradients, endows EXTRA with exact linear convergence to the optimal argument $\tbx^*$ under mild assumptions (\cite{Shi2014}).

If we recall the definitions of the local functions $f_{n}(\bbx_n)$ and the instantaneous local functions $f_{n,i}(\bbx_n)$ available at node $n$, the implementation of EXTRA requires that each node $n$ computes the {full} gradient of its local objective function $f_n$ at $\bbx_n^t$ as
\begin{equation}\label{gradient}
    \nabla f_{n}(\bbx_n^t) 
        = \frac{1}{q_n} \sum_{{i=1}}^{q_n}
               \nabla f_{n,{i}}(\bbx_n^t).
\end{equation}
This is computationally expensive when the number of instantaneous functions $q_n$ is large. To resolve this issue, local stochastic gradients can be substituted for the local objective functions gradients in \eqref{local_extra}. These stochastic gradients approximate the gradient $\nabla f_{n}(\bbx_n)$ of node $n$ by randomly choosing one of the instantaneous functions gradients $\nabla f_{n,i}(\bbx_n)$. If we let $i_n^t\in\{1,\ldots q_n\}$ denote a function index that we choose at time $t$ at node $n$ uniformly at random and independently of the history of the process, then the stochastic gradient is defined as 
\begin{equation}\label{sto_gradient}
    \hbs_n(\bbx_n^t):= \nabla f_{n,i_n^{t}}(\bbx_n^t).
\end{equation}
We can then write a stochastic version of EXTRA by replacing $\nabla f_n(\bbx_n^{t})$ by $\hbs_n(\bbx_n^t)$ and $\nabla f_n(\bbx_n^{t-1})$ by $\hbs_n(\bbx_n^{t-1})$. Such algorithm would have a small computational cost per iteration and, presumably, converge to the optimal argument $\tbx^*$. Here however, we want to design an algorithm with linear convergence rate, and stochastic descent algorithms achieve sublinear rates because of the difference between the stochastic and deterministic descent directions. 

%
\begin{figure} \centering
\ifcompiletikz \input{average_gradient_update.tex} 
\else \includegraphics[width=\linewidth]{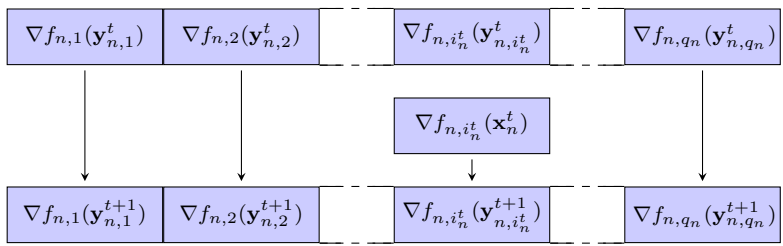}\fi 
\caption{Stochastic averaging gradient table at node $n$. At each iteration $t$ a random local instantaneous gradient $\nabla f_{n,i_n^{t}}(\bby_{n,i_{n}^t}^{t})$ is updated by $\nabla f_{n,i_n^{t}}(\bbx_n^{t})$. The rest of the local instantaneous gradients remain unchanged, i.e., $\nabla f_{n,i}(\bby_{n,i}^{t+1})=\nabla f_{n,i}(\bby_{n,i}^{t})$ for $i\neq i_n^t$. This list is used to compute the stochastic averaging gradient in \eqref{local_stochastic_gradient}.}
\label{fig_sag} \end{figure}

%
To reduce this noise we propose the use of stochastic averaging gradients instead (\cite{defazio2014saga}). The idea is to maintain a list of gradients of all instantaneous functions in which one randomly chosen element is replaced at each iteration and to use an average of the elements of this list for gradient approximation; see Figure \ref{fig_sag}. Formally, define the variable $\bby_{n,i}\in\reals^{p}$ to represent the iterate value the last time that the instantaneous gradient of function $f_{n,i}$ was evaluated. If we let $i_n^t\in\{1,\ldots, q_n\}$ denote the function index chosen at time $t$ at node $n$, as we did in \eqref{sto_gradient}, the variables $\bby_{n,i}$ are updated recursively as 
\begin{equation}\label{eqn_sag_update}
  \bby_{n,i}^{t+1} = \bbx_{n}^t,     \quad \text{if\ } i =    i_{n}^t, \qquad
  \bby_{n,i}^{t+1} = \bby_{n,i}^{t}, \quad \text{if\ } i \neq i_{n}^t.
\end{equation}
With these definitions in hand we can define the stochastic averaging gradient at node $n$ as
\begin{align}\label{local_stochastic_gradient}
   \hbg_{n}^t :=  \nabla f_{n,i_n^{t}}(\bbx_n^t)
                 -\nabla f_{n,i_n^{t}}(\bby_{n,i_{n}^t}^{t})
                 +\frac{1}{q_n}\sum_{i=1}^{q_n}
                  \nabla f_{n,i}(\bby_{n,i}^{t}).
\end{align}
Observe that to implement \eqref{local_stochastic_gradient} the gradients $\nabla f_{n,i}(\bby_{n,i}^{t})$ are stored in the local gradient table shown in Figure \ref{fig_sag}.

The DSA algorithm is a variation of EXTRA that substitutes the local gradients $\nabla f_n(\bbx_n^t)$ in \eqref{local_extra} for the local stochastic average gradients $\hbg_{n}^t$ in \eqref{local_stochastic_gradient},
\begin{equation}\label{DSA_local}
   \bbx_n^{t+1} =  \bbx_n^{t} 
                  +\sum_{m=1}^N w_{nm} \bbx_{m}^t 
                  -\sum_{m=1}^N \tdw_{nm} \bbx_{m}^{t-1}
                  - \alpha\left[\hbg_{n}^t-\hbg_{n}^{t-1}\right].
\end{equation}
The DSA initial update is given by applying the same substitution for the update of DGD in \eqref{gd_iteration} as
\begin{equation}\label{initial_DSA_local}
    \bbx_n^{1}=\sum_{m=1}^N w_{nm} \bbx_{m}^0 - \alpha\ \hbg_{n}^0.
\end{equation}
DSA is summarized in Algorithm \ref{algo_DSA} for $t\geq1$. The DSA update in \eqref{DSA_local} is implemented in Step 9. This step requires access to the local iterates $\bbx_m^t$ of  neighboring nodes $m\in \ccalN_{n}$ which are collected in Step 2. Furthermore, implementation of the DSA update also requires access to the stochastic averaging gradients $\hbg_n^{t-1}$ and $\hbg_n^t$. The latter is computed in Step 4 and the former is computed and stored at the same step in the previous iteration. The computation of the stochastic averaging gradients requires the selection of the index $i_{n}^t$. This index is chosen uniformly at random in Step 3. Determination of stochastic averaging gradients also necessitates access and maintenance of the gradients table in Figure \ref{fig_sag}. The $i_n^t$ element of this table is updated in Step 5 by replacing $\nabla f_{n,i_n^t}(\bby_{n,i_{n}^t}^t)$ with $\nabla f_{n,i_n^t}(\bbx_n^t)$, while the other vectors remain unchanged. To implement the first DSA iteration at time $t=0$ we have to perform the update in \eqref{initial_DSA_local} instead of the update in \eqref{DSA_local} as in Step 7. Furhter observe that the auxiliary variables $\bby_{n,i}^0$ are initialized to the initial iterate $\bbx_n^0$. This implies that the initial values of the stored gradients are $\nabla f_{n,i}(\bby_{n,i}^0) = \nabla f_{n,i}(\bbx_n^0)$ -- with a consequently relatively large initialization cost.

\begin{algorithm}[tb]
\caption{DSA algorithm at node $n$}\label{algo_DSA} 
\begin{algorithmic}[1] 
\small{\REQUIRE Vectors $\bbx_{n}^0$. Gradient table initialized with instantaneous gradients $\nabla f_{n,i}(\bby_{n,i}^0)$ with $\bby_{n,i}^0 = \bbx_n^0$. 
\FOR {$t=0,1,2,\ldots$}
   \STATE Exchange variable $\bbx_{n}^t$ with neighboring nodes $m\in \ccalN_n$.
   \STATE Choose $i_{n}^t$ uniformly at random from the set $\{1,\dots,q_{n}\}$.
   \STATE Compute and store stochastic averaging gradient as per \eqref{local_stochastic_gradient}:\\
          $ \displaystyle{\hbg_{n}^t
                 =  \nabla f_{n,i_n^{t}}(\bbx_{n}^t)
                  -\nabla f_{n,i_n^{t}}(\bby_{n,i_{n}^t}^{t}) 
                  +\frac{1}{q_n} \sum_{i=1}^{q_n} 
                   \nabla f_{n,i}(\bby_{n,i}^{t}) .}$
   \STATE Store $\nabla f_{n,i_{n}^t}(\bby_{n,i_{n}^t}^{t+1}) 
          = \nabla f_{n,i_{n}^t}(\bbx_{n}^t)$ in $i_{n}^t$ .
          gradient table position. %
  \IF{$t=0$}
  \vspace{-1.5mm}
     \STATE     Update variable $\bbx_{n}^t$ as per \eqref{initial_DSA_local}: 
     $\quad \displaystyle{\bbx_n^{t+1}=\sum_{n=1}^{N}w_{nm}\bbx_n^{t+1}- \alpha \hbg_{n}^t.}$
       \vspace{-1.5mm}
  \ELSE
    \vspace{-1.5mm}
     \STATE    Update variable $\bbx_{n}^t$ as per \eqref{DSA_local}:
          $\quad \displaystyle{\bbx_n^{t+1} 
              = \bbx_n^{t} + \sum_{n=1}^{N}w_{nm}\bbx_n^{t} 
                           - \sum_{n=1}^{N}\tilde{w}_{nm}\bbx_n^{t-1}
                           - \alpha\left[\hbg_{n}^t-\hbg_{n}^{t-1}\right].}$
                             \vspace{-1.5mm}
  \ENDIF
\ENDFOR}
\end{algorithmic}\end{algorithm}

%

We point out that the weights $w_{nm}$ and $\tdw_{nm}$ can't be arbitrary. If we define  weight matrices $\bbW$ and $\tbW$ with elements $w_{nm}$ and $\tdw_{nm}$, respectively, they have to satisfy conditions that we state as an assumption for future reference.
\\

%
\begin{assumption}\label{ass_wight_matrix_conditions}
The weight matrices $\bbW$ and $\tbW$ must satisfy the following properties
\begin{mylist}
\item[(a)] Both are symmetric, $\bbW=\bbW^T$ and $\tbW=\tbW^T$. \label{symmetry}
\item[(b)] The null space of $\bbI-\tbW$ includes the span of $\bbone$, i.e., $\text{null}(\bbI-\tbW)\supseteq \text{span}(\bbone)$, the null space of $\bbI-\bbW$ is the span of $\bbone$, i.e., $\text{null}(\bbI-\tbW)= \text{span}(\bbone)$, and the null space of the difference $\tbW-\bbW$ is the span of $\bbone$, i.e., $\text{null}(\tbW-\bbW)=\text{span}(\bbone)$. \label{spectral1}
\item[(c)] They satisfy the spectral ordering $ \bbW \preceq  \tbW \preceq ({\bbI+\bbW})/{2}$ and the matrix $\tbW$ is positive definite $\bb0\prec\tbW$.  \label{spectral2}
\end{mylist} \end{assumption}

%
Requiring the matrix $\bbW$ to be symmetric and with specific null space properties is necessary to let all agents converge to the same optimal variable. Analogous properties are necessary in DGD and are not difficult to satisfy. The condition on spectral ordering is specific to EXTRA but is not difficult to satisfy either. {E.g., if we have a matrix $\bbW$ that satisfies all the conditions in Assumption \ref{ass_wight_matrix_conditions}, the weight matrix $\tbW=(\bbI+\bbW)/2$ makes Assumption \ref{ass_wight_matrix_conditions} valid.} 

We also point that, as written in \eqref{local_stochastic_gradient}, computation of local stochastic averaging gradients $\hbg_n^t$ is costly because it requires evaluation of the sum $\sum_{i=1}^{q_n}\nabla f_{n,i}(\bby_{n,i}^{t})$ at each iteration. This cost can be avoided by updating the sum at each iteration with the recursive formula 
\begin{align}\label{local_stochastic_gradient_recursive}
 \sum_{i=1}^{q_n}\nabla f_{n,i}(\bby_{n,i}^{t}) 
 	   \ =\  \sum_{i=1}^{q_n}\nabla f_{n,i}(\bby_{n,i}^{t-1}) 
           + \nabla f_{n,i_n^{t-1}}(\bbx_n^{t-1})
                        -\nabla f_{n,i_n^{t-1}}(\bby_{n,i_{n}^{t-1}}^{t-1}) .
\end{align}
Important properties and interpretations of EXTRA and DSA are presented in the following sections after a pertinent remark.

%
\begin{remark}\label{rmk_unbiased_estimation}\normalfont The local stochastic averaging gradients in \eqref{local_stochastic_gradient} are unbiased estimates of the local gradients $\nabla f_{n}(\bbx_n^t)$. Indeed, if we let $\ccalF_t$ measure the history of the system up until time $t$ we have that the sum in \eqref{local_stochastic_gradient} is deterministic given this sigma-algebra. Thus, the conditional expectation of the stochastic averaging gradient is, 
\begin{align}\label{expected_value_simplification_intermediate}
   \E{\hbg_{n}^t\given\ccalF^t} 
       = \mbE\Big[\nabla f_{n,i_n^{t}}(\bbx_n^t)\given\ccalF^t\Big]
                 -\mbE\Big[\nabla f_{n,i_n^{t}}(\bby_{n,i_{n}^t}^{t})\given\ccalF^t\Big]
                 +\frac{1}{q_n}\sum_{i=1}^{q_n}
                  \nabla f_{n,i}(\bby_{n,i}^{t}).
\end{align}
With the index $i_n^t$ chosen equiprobably from the set $\{1,\ldots, q_n\}$, the expectation of the second term in \eqref{expected_value_simplification_intermediate} is the same as the sum in the last term -- each of the indexes is chosen with probability $1/q_n$. Therefore, these two terms cancel out each other and, since the expectation of the first term in \eqref{expected_value_simplification_intermediate} is simply $\mbE\big[\nabla f_{n,i_n^{t}}(\bbx_n^t)\given\ccalF^t\big] = (1/q_n) \sum_{i=1}^{q_n}\nabla f_{n,i}(\bbx_n^t) =\nabla f_{n}(\bbx_n^t)$, we can simplify \eqref{expected_value_simplification_intermediate} to
\begin{align}\label{expected_value_simplification}
   \E{\hbg_{n}^t\given\ccalF^t} = \nabla f_{n}(\bbx_n^t).
\end{align}
The expression in \eqref{expected_value_simplification} means, by definition, that $\hbg_{n}^t$ is an unbiased estimate of $\nabla f_{n}(\bbx_n^t)$ when the history $\ccalF^t$ is given.
\end{remark}

%% file: average_gradient_update.tex
\def \thisplotscale {0.74}
\def \unit {\thisplotscale cm}

%
\tikzstyle{block}         = [ draw,
                              rectangle,
                              minimum height = 1.0*\unit,
                              minimum width  = 2.5*\unit,
                              text width     = 2.5*\unit,
                              text badly centered,
                              fill = blue!20, 
                              font = \footnotesize, 
                              anchor = west]
\tikzstyle{bold block}    = [ block,
                              fill = blue!40]
\tikzstyle{light block}   = [ block,
                              fill = blue!10]
\tikzstyle{connector}     = [ draw, 
                              -stealth, 
                              shorten >=2,
                              shorten <=2,]
\tikzstyle{dot dot dot}   = [ draw, 
                              dashed]

%
{\small
 \begin{tikzpicture}[x = 0.95*\unit, y = 0.8*\unit]
  

\path (0,0)                 node [block] (S1) {$\nabla f_{n,1}(\bby_{n,1}^t)$}; 
\path (S1.east)             node [block] (S2) {$\nabla f_{n,2}(\bby_{n,2}^t)$}; 
\path (S2.east) ++ (1.4, 0) node [block] (Sr) {$\nabla f_{n,i_n^t}(\bby_{n,i_n^t}^t)$}; 
\path (Sr.east) ++ (1.4, 0) node [block] (Sq) {$\nabla f_{n,q_n}(\bby_{n,q_n}^t)$}; 
\path (Sr.west) ++ (0,-2.0) node [block] (Ir) {$\nabla f_{n,i_n^t}(\bbx_n^t)$}; 
\path (S1.west) ++ (0,-4.0) node [block] (D1) {$\nabla f_{n,1}(\bby_{n,1}^{t+1})$}; 
\path (D1.east)             node [block] (D2) {$\nabla f_{n,2}(\bby_{n,2}^{t+1})$}; 
\path (D2.east) ++ (1.4, 0) node [block] (Dr) {$\nabla f_{n,i_n^t}(\bby_{n,i_n^{t}}^{t+1})$}; 
\path (Dr.east) ++ (1.4, 0) node [block] (Dq) {$\nabla f_{n,q_n}(\bby_{n,q_n}^{t+1})$}; 

\path[draw]        (S2.north east) -- ++ ( 0.4, 0);
\path[draw]        (Sr.north west) -- ++ (-0.4, 0);
\path[dot dot dot] (S2.north east) -- (Sr.north west);
\path[draw]        (Sr.north east) -- ++ ( 0.4, 0);
\path[draw]        (Sq.north west) -- ++ (-0.4, 0);
\path[dot dot dot] (Sr.north east) -- (Sq.north west);
\path[draw]        (S2.south east) -- ++ ( 0.4, 0);
\path[draw]        (Sr.south west) -- ++ (-0.4, 0);
\path[dot dot dot] (S2.south east) -- (Sr.south west);
\path[draw]        (Sr.south east) -- ++ ( 0.4, 0);
\path[draw]        (Sq.south west) -- ++ (-0.4, 0);
\path[dot dot dot] (Sr.south east) -- (Sq.south west);

\path[draw]        (D2.north east) -- ++ ( 0.4, 0);
\path[draw]        (Dr.north west) -- ++ (-0.4, 0);
\path[dot dot dot] (D2.north east) -- (Dr.north west);
\path[draw]        (Dr.north east) -- ++ ( 0.4, 0);
\path[draw]        (Dq.north west) -- ++ (-0.4, 0);
\path[dot dot dot] (Dr.north east) -- (Dq.north west);
\path[draw]        (D2.south east) -- ++ ( 0.4, 0);
\path[draw]        (Dr.south west) -- ++ (-0.4, 0);
\path[dot dot dot] (D2.south east) -- (Dr.south west);
\path[draw]        (Dr.south east) -- ++ ( 0.4, 0);
\path[draw]        (Dq.south west) -- ++ (-0.4, 0);
\path[dot dot dot] (Dr.south east) -- (Dq.south west);

\path[connector] (S1.south) -- (D1.north);
\path[connector] (S2.south) -- (D2.north);
\path[connector] (Ir.south) -- (Dr.north);
\path[connector] (Sq.south) -- (Dq.north);

\end{tikzpicture}}

%% file: saddle_point_interpretation.tex

%
\subsection{Limit points of DGD and EXTRA}\label{subset:lmit_point}

The derivation of EXTRA hinges on the observation that the optimal argument of \eqref{opt_problem} is not a fixed point of the DGD iteration in \eqref{gd_iteration} but is a fixed point of the iteration in \eqref{local_extra}. To explain this point define $\bbx := \left[\bbx_{1}; \dots ; \bbx_{N}\right]\in \reals^{Np}$ as a vector that concatenates the local iterates $\bbx_{n}$ and the aggregate function $f:\reals^{Np}\to \reals$ as the one that takes values $f(\bbx)=f(\bbx_1,\dots, \bbx_N):=\sum_{n=1}^N f_{n}(\bbx_n)$. Decentralized optimization entails the minimization of $f(\bbx)$ subject to the constraint that all local variables are equal,
\begin{alignat}{2}\label{agg_func}
   \bbx^* :=
    &\argmin &&\ f\left(\bbx\right)
                    = f(\bbx_1,\dots, \bbx_N)
                    = \sum_{n=1}^N f_{n}(\bbx_n), \nonumber\\
    &\st     &&\ \bbx_n = \bbx_m ,\quad  \text{for all\ } n,m.
\end{alignat}
The problems in \eqref{opt_problem} and \eqref{agg_func} are equivalent in the sense that the vector $\bbx^*\in \reals^{Np}$ is a solution of \eqref{agg_func} if it satisfies $\bbx_n^* = \tbx^*$ for all $n$, or, equivalently, if we can write $\bbx^*=[\tbx^*;\dots;\tbx^*]$. Regardless of interpretation, the Karush, Kuhn, Tucker (KKT) conditions of  \eqref{agg_func} dictate that that optimal argument $\bbx^*$ must sastisfy
\begin{equation}\label{opt_cond}
   \bbx^* \subset \text{span}(\bold{1}_{N}\otimes\bbI_p) , \qquad 
   (\bold{1}_{N}\otimes\bbI_p)^T\nabla f(\bbx^*) = \bbzero.
\end{equation}
The first condition in \eqref{opt_cond} requires that all the local variables $\bbx_n^*$ be equal, while the second condition requires the sum of local gradients to vanish at the optimal point. This latter condition is not the same as $\nabla f(\bbx)=\bbzero$. If we observe that the gradient $\nabla f(\bbx^t)$ of the aggregate function can be written as $\nabla f(\bbx)=[\nabla f_1(\bbx_1);\dots;\nabla f_N(\bbx_N)]\in \reals^{Np}$, the condition $\nabla f(\bbx)=\bbzero$ implies that all the local gradients are null, i.e., that $\nabla f_{n}(\bbx_n)=\bb0$ for all $n$. This is stronger than having their sum being null as required by \eqref{opt_cond}.

Define now 
the extended weight matrices as the Kronecker products $\bbZ:= \bbW \otimes \bbI\in\reals^{Np\times Np}$ and $\tbZ:= \tbW \otimes \bbI\in\reals^{Np\times Np}$. Note that the required conditions for the weight matrices $\bbW$ and $\tbW$ in Assumption \ref{ass_wight_matrix_conditions} enforce some conditions on the extended weight matrices $\bbZ$ and $\tbZ$. Based on Assumption \ref{ass_wight_matrix_conditions}(a), the matrices $\bbZ$ and $\tbZ$ are also symmetric, i.e., $\bbZ=\bbZ^T$ and $\tbZ=\tbZ^T$. Conditions in Assumption \ref{ass_wight_matrix_conditions}(b) imply that null$\{\tbZ-\bbZ\}=$ span$\{\bold{1}\otimes \bbI\}$, null$\{\bbI-\bbZ\}=$ span$\{\bold{1}\otimes \bbI\}$, and null$\{\bbI-\tbZ\}\supseteq$ span$\{\bold{1}\otimes \bbI\}$. Lastly, the spectral properties of matrices $\bbW$ and $\tbW$ in Assumption \ref{ass_wight_matrix_conditions}(c) yield that matrix $\tbZ$ is positive definite and the expression $\bbZ \preceq  \tbZ \preceq ({\bbI+\bbZ})/{2}$ holds.

According to the definition of extended weight matrix $\bbZ$, the DGD iteration in \eqref{gd_iteration} is equivalent to
\begin{equation}\label{new_dgd_update}
 \bbx^{t+1} = \bbZ \bbx^t - \alpha\nabla f(\bbx^t),
\end{equation}
where, according to \eqref{agg_func}, the gradient $\nabla f(\bbx^t)$ of the aggregate function can be written as $\nabla f(\bbx^t)=[\nabla f_1(\bbx_1^t);\dots;\nabla f_N(\bbx_N^t)]\in \reals^{Np}$. Likewise, 
the EXTRA iteration in \eqref{local_extra} can be written as
\begin{equation}\label{extra_update}
    \bbx^{t+1} = (\bbI+\bbZ)\bbx^{t} -\tbZ\bbx^{t-1} 
                       - \alpha\left[{\nabla} f(\bbx^{t})-{\nabla} f(\bbx^{t-1})\right].
\end{equation}
The fundamental difference between DGD and EXTRA is that a fixed point of \eqref{new_dgd_update} does not necessarily satisfy \eqref{opt_cond}, whereas the fixed points of \eqref{extra_update} are guaranteed to do so. Indeed, taking limits in \eqref{new_dgd_update} we see that the fixed points $\bbx^\infty$ of DGD must satisfy
\begin{equation}\label{eqn_dgd_limit_point}
   (\bbI-\bbZ) \bbx^{\infty} + \alpha\nabla f(\bbx^\infty) =\bbzero,
\end{equation}
which is incompatible with \eqref{opt_cond} except in peculiar circumstances -- such as, e.g., when all local functions have the same minimum. The limit points of EXTRA, however, satisfy the relationship
\begin{equation}\label{extra_limit_point}
    \bbx^{\infty}-\bbx^{\infty}
        = (\bbZ- \tbZ)\bbx^{\infty} 
          - \alpha[{\nabla} f(\bbx^{\infty})-{\nabla} f(\bbx^{\infty})].
\end{equation}
Canceling out the variables on the left hand side and the gradients in the right hand side it follows that $(\bbZ- \tbZ)\bbx^{\infty}=\bbzero$. Since the null space of of $\bbZ- \tbZ$ is $\text{null}(\bbZ- \tbZ)=\bold{1}_{N}\otimes\bbI_p$ by assumption, we must have $\bbx^\infty \subset \text{span}(\bold{1}_{N}\otimes\bbI_p)$. This is the first condition in \eqref{opt_cond}. For the second condition in \eqref{opt_cond} sum the updates in \eqref{extra_update} recursively and use the telescopic nature of the sum to write
\begin{equation}\label{extra_new_relation}
   \bbx^{t+1}=\tbZ\bbx^{t} - \alpha{\nabla} f(\bbx^{t}) -\sum_{s=0}^{t}(\tbZ-\bbZ)\bbx^s.
\end{equation}
Substituting the limit point in \eqref{extra_new_relation} and reordering terms, we see that $\bbx^{\infty}$ must satisfy
\begin{equation}\label{extra_limit_point_second_relation}
   \alpha{\nabla} f(\bbx^{\infty}) 
      = (\bbI-\tbZ) \bbx^{\infty} 
          - \sum_{s=0}^{\infty}(\tbZ-\bbZ)\bbx^s.
\end{equation}
In \eqref{extra_limit_point_second_relation} we have that $(\bbI-\tbZ) \bbx^{\infty} =\bbzero$ because the null space of $(\bbI-\tbZ)$ is $\text{null}(\bbZ- \tbZ)=\bold{1}_{N}\otimes\bbI_p$ by assumption and $\bbx^\infty \subset \text{span}(\bold{1}_{N}\otimes\bbI_p)$ as already shown. Implementing this simplification and considering the multiplication of the resulting equality by $(\bold{1}_{N}\otimes\bbI_p)^T$ we obtain
\begin{equation}\label{extra_limit_point2}
    (\bold{1}_{N}\otimes\bbI_p)^T \alpha{\nabla} f(\bbx^{\infty}) 
        = -\sum_{s=0}^{\infty}(\bold{1}_{N}\otimes\bbI_p)^T(\bbZ-\tbZ)\bbx^s.
\end{equation}
In \eqref{extra_limit_point2}, the terms $(\bold{1}_{N}\otimes\bbI_p)^T(\bbZ-\tbZ)=\bbzero$ because the matrices $\bbZ$ and $\tbZ$ are symmetric and $(\bold{1}_{N}\otimes\bbI_p)$ is in the null space of the difference $\bbZ-\tbZ$. This implies that $(\bold{1}_{N}\otimes\bbI_p)^T \alpha{\nabla} f(\bbx^{\infty})$, which is the second condition in \eqref{agg_func}. Therefore, given the assumption that the sequence of EXTRA iterates $\bbx^{t}$ has a limit point $\bbx^{\infty}$ it follows that this limit point satisfies both conditions in \eqref{opt_cond} and for this reason exact convergence with constant stepsize is achievable for EXTRA.

%
\subsection{Stochastic saddle point method interpretation of DSA}\label{subset:saddle}

The convergence proofs of DSA build on a reinterpretation of EXTRA as a saddle point method. To introduce this primal-dual interpretation consider the update in \eqref{extra_new_relation} and define the sequence of vectors $\bbv^t=\sum_{s=0}^{t}(\tbZ-\bbZ)^{1/2}\bbx^s$. The vector $\bbv^t$ represents the accumulation of variable dissimilarities in different nodes over time. Considering this definition of $\bbv^t$ we can rewrite \eqref{extra_new_relation} as
\begin{equation}\label{extra_primal_update}
\bbx^{t+1}=\bbx^{t} -\alpha\! \left[{\nabla} f(\bbx^{t})+\frac{1}{\alpha}(\bbI\!-\!\tbZ)\bbx^t  + \frac{1}{\alpha}(\tbZ\!-\!\bbZ)^{1/2}\bbv^t\right].
\end{equation}
Furthermore, based on the definition of the sequence  $\bbv^t=\sum_{s=0}^{t}(\tbZ-\bbZ)^{1/2}\bbx^s$ we can write the recursive expression
\begin{equation}\label{extra_dual_update}
\bbv^{t+1}=\bbv^{t} + \alpha\left[ \frac{1}{\alpha}(\tbZ-\bbZ)^{1/2}\bbx^{t+1}\right].
\end{equation}
Consider $\bbx$ as a primal variable and $\bbv$ as a dual variable. Then, the updates in \eqref{extra_primal_update} and \eqref{extra_dual_update} are equivalent to the updates of a saddle point method with stepsize $\alpha$ that solves for the critical points of the augmented Lagrangian 
\begin{equation}\label{lagrangian}
\ccalL(\bbx,\bbv)= f(\bbx)+\frac{1}{\alpha}\bbv^T(\tbZ-\bbZ)^{1/2}\bbx+\frac{1}{2\alpha}\bbx^T(\bbI-\tbZ)\bbx.
\end{equation}
In the Lagrangian in \eqref{lagrangian} the factor $(1/\alpha)\bbv^T(\tbZ-\bbZ)^{1/2}\bbx$ stems from the linear constraint $(\tbZ-\bbZ)^{1/2}\bbx=\bbzero$ and the quadratic term $({1}/{2\alpha})\bbx^T(\bbI-\tbZ)\bbx$ is the augmented term added to the Lagrangian. Therefore, the optimization problem whose augmented Lagrangian is the one given in \eqref{lagrangian} is
\begin{equation}\label{constrained_opt_problem}
\bbx^* =\argmin_{\bbx} \ f(\bbx)\qquad \text{s.t.} \ \frac{1}{\alpha}(\tbZ-\bbZ)^{1/2}\bbx=\bb0.
\end{equation}
Observing that the null space of $(\tbZ-\bbZ)^{1/2}$ is $\text{null}((\tbZ-\bbZ)^{1/2})=\text{null}{(\tbZ-\bbZ)}=\text{span}\{\bold{1}_{N}\otimes\bbI_p\}$, the constraint in \eqref{constrained_opt_problem} is equivalent to the consensus constraint $\bbx_n = \bbx_m$ for all $n,m$ that appears in \eqref{agg_func}. This means that \eqref{constrained_opt_problem} is equivalent to \eqref{agg_func}, which, as already argued, is equivalent to the original problem in \eqref{opt_problem}. Hence, EXTRA is a saddle point method that solves \eqref{constrained_opt_problem} which, because of their equivalence, is tantamount to solving \eqref{opt_problem}. Considering that saddle point methods converge linearly, it follows that the same is true of EXTRA.

That EXTRA is a saddle point method provides a simple explanation of its convergence properties. For the purposes of this paper, however, the important fact is that if EXTRA is a saddle point method, DSA is a stochastic saddle point method. To write DSA in this form define $\hbg^t:=[\hbg_{1}^t;\dots;\hbg_{N}^t]\in \reals^{Np}$ as the vector that concatenates all the local stochastic averaging gradients at step $t$. Then, the DSA  update in \eqref{DSA_local} can be written as
\begin{align}\label{DSA_update}
   \bbx^{t+1}=(\bbI+\bbZ)\bbx^{t} -\tbZ\bbx^{t-1}- \alpha\left[\hbg^{t}-\hbg^{t-1}\right].
\end{align}
Comparing \eqref{extra_update} and \eqref{DSA_update} we see that they differ in the latter  using stochastic averaging gradients $\hbg^t$ in lieu of the full gradients $\nabla f(\bbx^t)$. Therefore, DSA is a stochastic saddle point method in which the primal variables are updated as
\begin{equation}\label{DSA_primal_update}
   \bbx^{t+1}=\bbx^{t} -\alpha \hbg^{t}-(\bbI-\tbZ)\bbx^t  -(\tbZ-\bbZ)^{1/2}\bbv^t,
\end{equation}
and the dual variables $\bbv^t$ are updated as 
\begin{equation}\label{DSA_dual_update}
\bbv^{t+1}=\bbv^{t} + (\tbZ-\bbZ)^{1/2}\bbx^{t+1}.
\end{equation}
Notice that the initial primal variable $\bbx^0$ is an arbitrary vector in $\reals^{Np}$, while according to the definition $\bbv^{t}=\sum_{s=0}^t(\tbZ-\bbZ)^{1/2}\bbx^{s}$. We then need to set the initial multiplier to $\bbv^{0}=(\tbZ-\bbZ)^{1/2}\bbx^{0}$. This is not a problem in practice because \eqref{DSA_primal_update} and \eqref{DSA_dual_update} are not used for implementation. In our converge analysis we utilize the (equivalent) stochastic saddle point expressions for DSA shown in \eqref{DSA_primal_update} and \eqref{DSA_dual_update}. The expression in \eqref{DSA_local} is used for implementation because it avoids exchanging dual variables -- as well as the initialization problem. The convergence analysis is presented in the following section.

%% file: Convergence2.tex
\section{Convergence analysis}\label{sec_convergence}

Our goal here is to show that as time progresses the sequence of iterates $\bbx^t$ approaches the optimal argument $\bbx^*$. To do so, in addition to the conditions on the weight matrices $\bbW$ and $\tbW$ in Assumption \ref{ass_wight_matrix_conditions}, we assume the instantaneous local functions $f_{n,i}$ have specific properties that we state next.



%
\begin{assumption}\label{ass_intantaneous_hessian_bounds}\normalfont  The instantaneous local functions $f_{n,i}(\bbx_n)$ are differentiable and strongly convex with parameter $\mu$.
\end{assumption}

%
\begin{assumption}\normalfont\label{ass_bounded_stochastic_gradient_norm} 
The gradient of instantaneous local functions $\nabla f_{n,i}$ are Lipschitz continuous with parameter $L$. I.e., for all $n\in\{1,\dots,N\}$ and $i\in\{1,\dots,q_n\}$ we can write
\begin{equation}
\left\|\nabla f_{n,i}(\bba)-\nabla f_{n,i}(\bbb)\right\| \leq L\ \| \bba-\bbb \|\quad\!  \bba,\bbb \in \reals^p.
\end{equation}
\end{assumption}

The condition imposed by Assumption \ref{ass_intantaneous_hessian_bounds} implies that the local functions $f_{n}(\bbx_{n})$ and the global cost function $f(\bbx)=\sum_{n=1}^N f_{n}(\bbx_{n})$ are also strongly convex with parameter $\mu$. Likewise, Lipschitz continuity of the local instantaneous gradients considered in Assumption \ref{ass_bounded_stochastic_gradient_norm} enforces Lipschitz continuity of gradients of the local functions $\nabla f_n(\bbx_n)$ and the aggregate function $\nabla f(\bbx)$ -- see, e.g., (Lemma 1 of \cite{mokhtari2015network}).

\subsection{Preliminaries}

In this section we study some basic properties of the sequences of primal and dual variables generated by the DSA algorithm. In the following lemma, we study the relation of iterates $\bbx^t$ and $\bbv^t$ with the optimal primal  $\bbx^*$ and dual $\bbv^*$ arguments.

\begin{lemma}\label{lemma_important}
Consider the DSA algorithm as defined in \eqref{eqn_sag_update}-\eqref{initial_DSA_local} and recall the updates of the primal $\bbx^t$ and dual $\bbv^t$ variables in \eqref{DSA_primal_update} and \eqref{DSA_dual_update}, respectively. Further, define the positive semidefinite matrix $\bbU:=(\tbZ-\bbZ)^{1/2}$. If Assumption \ref{ass_wight_matrix_conditions} holds true, then the sequence of primal $\bbx^t$ and dual $\bbv^t$ variables satisfy 
\begin{align}\label{claim_10}
\alpha \left[\hbg^{t} -\nabla f(\bbx^*)  \right] =(\bbI+\bbZ-2\tbZ)(\bbx^*-\bbx^{t+1})+\tbZ(\bbx^{t}-\bbx^{t+1}) -\bbU(\bbv^{t+1}-\bbv^*).
\end{align}
\end{lemma}

\begin{proof}
Considering the update rule for the dual variable in \eqref{DSA_dual_update} and the definition $\bbU=(\tbZ-\bbZ)^{1/2}$, we can substitute $\bbU\bbv^t$ in \eqref{DSA_primal_update} by $\bbU\bbv^{t+1}-\bbU^2\bbx^{t+1}$. Applying this substitution into the DSA primal update in \eqref{DSA_primal_update} yields
\begin{equation}\label{combine_20}
\alpha \hbg^{t}
= -(\bbI+\bbZ-\tbZ)\bbx^{t+1}+\tbZ\bbx^{t}-\bbU\bbv^{t+1}.
\end{equation}
By adding and subtracting $\tbZ\bbx^{t+1}$ into the right hand side of \eqref{combine_20} and considering the fact that $(\bbI+\bbZ-2\tbZ)\bbx^{*}=\bb0$ we obtain 
\begin{equation}\label{combine_30}
\alpha \hbg^{t}
= (\bbI+\bbZ-2\tbZ)(\bbx^*-\bbx^{t+1})+\tbZ(\bbx^{t}-\bbx^{t+1})-\bbU\bbv^{t+1}.
\end{equation}
One of the KKT conditions of problem \eqref{constrained_opt_problem} follows that the optimal variables $\bbx^*$ and $\bbv^*$ satisfy $\alpha\nabla f(\bbx^*)+\bbU\bbv^*=\bb0$ or equivalently $-\alpha \nabla f(\bbx^*)=\bbU\bbv^{*}$. Adding this equality to both sides of \eqref{combine_30} follows the claim in \eqref{claim_10}. 
\end{proof}


In the subsequent analyses of convergence of DSA, we need an upper bound for the expected value of squared difference between the stochastic averaging gradient $\hbg^t$ and the gradient of optimal argument $\nabla f(\bbx^*)$ given the observation until step $t$, i.e. $\E{\left\|   \hbg^{t}    -    \nabla f(\bbx^*) \right\|^2  \mid \ccalF^t}$. To establish this upper bound first we define the sequence $p^{t}\in \reals$ as
\begin{align}\label{p_definition}
p^{t}:= \sum_{n=1}^N
\bigg[ 
\frac{1}{q_n}\sum_{i=1}^{q_n}\left( f_{n,i}(\bby_{n,i}^{t}) - f_{n,i}(\tbx^*) -\nabla f_{n,i}(\tbx^*)^T (\bby_{n,i}^{t}-\tbx^*)\right)
 \bigg].
\end{align}
Notice that based on strong convexity of local instantaneous functions $f_{n,i}$, each term $f_{n,i}(\bby_{n,i}^{t}) - f_{n,i}(\tbx^*) 
-\nabla f_{n,i}(\tbx^*)^T (\bby_{n,i}^{t}-\tbx^*)$ is positive and as a result the sequence $p^t$ defined in \eqref{p_definition} is always positive. In the following lemma, we use the result in Lemma \ref{lemma_important} to guarantee an upper bound for the expectation $\E{\|  \hbg^{t} -\nabla f(\bbx^*) \|^2  \mid \ccalF^t}$ in terms of $p^t$ and the optimality gap $f(\bbx^{t}) - f(\bbx^*)  -\nabla f(\bbx^*)^T (\bbx^{t}-\bbx^*)$.

\begin{lemma}\label{lemma_gradient_noise_bound}
Consider the DSA algorithm in \eqref{eqn_sag_update}-\eqref{initial_DSA_local} and the definition of sequence $p^t$ in \eqref{p_definition}. If Assumptions \ref{ass_wight_matrix_conditions}-\ref{ass_bounded_stochastic_gradient_norm} hold true, then the squared norm of the difference between stochastic averaging gradient  $\hbg^{t}$ and the optimal gradient $\nabla f(\bbx^*)$ in expectation is bounded above by
\begin{align}\label{claim_gradient_noise}
\E{\left\|   \hbg^{t}    -    \nabla f(\bbx^*) \right\|^2  
\!\mid\! \ccalF^t}	 \leq 
			4Lp^t
+2\left(2 L -\mu\right)  \left(f(\bbx^{t}) - f(\bbx^*)  -\nabla f(\bbx^*)^T (\bbx^{t}-\bbx^*)\right).
\end{align}
\end{lemma}

\begin{proof}
See Appendix \ref{apx_lemma_grad_noise_bound}.
\end{proof}

Observe that as the sequence of iterates $\bbx^t$ approaches the optimal argument $\bbx^*$, all the local auxiliary variables $\bby_{n,i}^t$ converge to $\tbx^*$ which follows convergence of $p^t$ to null. This observation in association with the result in \eqref{claim_gradient_noise} implies that the expected value of the difference between the stochastic averaging gradient $\hbg^{t}  $ and the optimal gradient $\nabla f(\bbx^*)$ vanishes as the sequence of iterates $\bbx^t$ approaches the optimal argument $\bbx^*$.


\subsection{Convergence}

In this section we establish linear convergence of the sequence of iterates $\bbx^t$ generated by DSA to the optimal argument $\bbx^*$. To do so define $0<\gamma $ and $\Gamma<\infty$ as the smallest and largest eigenvalues of positive definite matrix matrix $\tbZ$, respectively. Likewise, define $\gamma'$ as the smallest non-zero eigenvalue of matrix $\tbZ-\bbZ$ and $\Gamma'$ as the largest eigenvalue of matrix $\tbZ-\bbZ$. Further, define vectors $\bbu^{t},\bbu^*\in \reals^{2Np}$ and matrix $\bbG\in \reals^{2Np\times 2Np}$ as
\begin{equation}\label{z_G_definitions}
\bbu^*:=
\begin{bmatrix}
\bbx^*\\ \bbv^*
\end{bmatrix}
, \quad 
\bbu^t:=
\begin{bmatrix}
\bbx^t\\ \bbv^t
\end{bmatrix}
,\quad
\bbG
=
\begin{bmatrix}
    \tbZ & \bb0 \\
    \bb0 & \bbI
\end{bmatrix}.
\end{equation}
Vector $\bbu^*\in \reals^{2Np}$ concatenates the optimal primal and dual variables and vector $\bbu^t\in \reals^{2Np}$ contains primal and dual iterates at step $t$. Matrix $\bbG\in \reals^{2Np\times 2Np}$ is a block diagonal positive definite matrix that we introduce since  instead of tracking the value of $\ell_2$ norm  $\|\bbu^{t}-\bbu^{*}\|^2_{2}$ we study the convergence properties of $\bbG$ weighted norm $\|\bbu^{t}-\bbu^{*}\|^2_{\bbG}$. Notice that the weighted norm $\|\bbu^{t}-\bbu^{*}\|^2_{\bbG}$ is equivalent to $(\bbu^{t}-\bbu^{*})^T\bbG(\bbu^{t}-\bbu^{*})$. Our goal is to show that the sequence $\|\bbu^{t}-\bbu^{*}\|^2_{\bbG}$ converges linearly to null. To do this we show linear convergence of a Lyapunov function of the sequence $\|\bbu^{t}-\bbu^{*}\|^2_{\bbG}$. The Lyapunov function is defined as $\|\bbu^{t}-\bbu^{*}\|^2_{\bbG}+c p^t$ where $c>0$ is a positive constant.

To prove linear convergence of the sequence $\|\bbu^{t}-\bbu^{*}\|^2_{\bbG}+c p^t$ we first show an upper bound for the expected error $\E{\|\bbu^{t+1}-\bbu^{*}\|^2_{\bbG}\mid \ccalF^t}$ in terms of $\|\bbu^{t}-\bbu^{*}\|^2_{\bbG}$ and some parameters that capture the optimality gap.

\begin{lemma}\label{lemma_lower_bound_for_z_decrement}
Consider the DSA algorithm as defined in \eqref{eqn_sag_update}-\eqref{initial_DSA_local}. Further recall the definitions of $p^t$ in \eqref{p_definition} and $\bbu^t$, $\bbu^*$, and $\bbG$ in \eqref{z_G_definitions}. If Assumptions \ref{ass_wight_matrix_conditions}-\ref{ass_bounded_stochastic_gradient_norm} hold true, then for any positive constants $\eta>0$ we can write 
\begin{align}\label{claim_z_decrement_prime}
	\E{\|\bbu^{t+1}-\bbu^*\|_\bbG^2\mid \ccalF^t}	&\leq
\|\bbu^t-\bbu^*\|_\bbG^2
-2\E{\|\bbx^{t+1}-\bbx^*\|_{\bbI+\bbZ-2\tbZ}^2\mid \ccalF^t}
	+\frac{\alpha4L}{\eta}p^t
\nonumber
\\
&\quad 	-\E{\|\bbx^{t+1}-\bbx^t\|_{\tbZ-2\alpha\eta\bbI}^2\mid \ccalF^t}
	-\E{\|\bbv^{t+1}-\bbv^t\|^2\mid \ccalF^t}\nonumber\\
&   \quad  
	 -\left(\frac{4\alpha\mu }{L} -   \frac{2\alpha( 2L-\mu)}{\eta}  \right)
\left(
 f(\bbx^{t}) - f(\bbx^*)  -\nabla f(\bbx^*)^T (\bbx^{t}-\bbx^*)\right).
\end{align}
\end{lemma}

\begin{proof}
See Appendix \ref{apx_lemma_lower_bound_for_z_decrement}.
\end{proof}


Lemma \ref{lemma_lower_bound_for_z_decrement} shows an upper bound for the squared norm $\|\bbu^{t+1}-\bbu^{*}\|^2_{\bbG}$ which is the first part of the Lyapunov function $\|\bbu^{t}-\bbu^{*}\|^2_{\bbG}+c p^t$ at step $t+1$. Likewise, we provide an upper bound for the second term of the Lyapunov function at time $t+1$ which is $p^{t+1}$ in terms of $p^{t}$ and some parameters that capture optimality gap. This bound is studied in the following lemma.

\begin{lemma}\label{lemma_p_decrement_2}
Consider the DSA algorithm as defined in \eqref{eqn_sag_update}-\eqref{initial_DSA_local} and the definition of $p^t$ in \eqref{p_definition}. Further, define $q_{\min}$ and $q_{\max}$ as the smallest and largest values for the number of instantaneous functions at a node, respectively. If Assumptions \ref{ass_wight_matrix_conditions}-\ref{ass_bounded_stochastic_gradient_norm} hold true, then for all $t>0$ the sequence $p^t$ satisfies
\begin{align}\label{lemma_p_decrement_claim}
\E{
p^{t+1}
 \mid \ccalF^t}
 \leq
 \left[1-\frac{1}{q_{\max}}\right]p^t+\frac{1}{q_{\min}}\left[ f(\bbx^t)-f(\bbx^*) -\nabla f(\bbx^*)^T (\bbx^{t}-\bbx^*) \right].
\end{align}
\end{lemma}

\begin{proof}
See Appendix \ref{apx_lemma_p_decrement_2}.
\end{proof}

Lemma \ref{lemma_p_decrement_2} provides an upper bound for $p^{t+1}$ in terms of its previous value $p^{t}$ and the optimality error $f(\bbx^t)-f(\bbx^*) -\nabla f(\bbx^*)^T (\bbx^{t}-\bbx^*) $. Combining the results in Lemmata \ref{lemma_lower_bound_for_z_decrement} and \ref{lemma_p_decrement_2} we can show that in expectation the Lyapunov function $\| \bbu^{t+1}-\bbu^*\|_\bbG^2+c\ p^{t+1}$ at step $t+1$ is strictly smaller than its previous value  $\| \bbu^{t}-\bbu^*\|_\bbG^2+c\ p^t$ at step $t$.

\begin{theorem}\label{lin_convg_thm}
Consider the DSA algorithm as defined in \eqref{eqn_sag_update}-\eqref{initial_DSA_local}. Further recall the definition of the sequence $p^t$ in \eqref{p_definition}. Define $\eta$ as an arbitrary positive constant chosen from the interval 
\begin{equation}\label{eta_condition}
\eta\in\left(
 \frac{L^2q_{\max}}{\mu q_{\min} }+\frac{L^2}{\mu}-L\ , \ \infty
\right).
\end{equation}
If Assumptions \ref{ass_wight_matrix_conditions}-\ref{ass_bounded_stochastic_gradient_norm} hold true and the stepsize $\alpha$ is chosen from the interval $\alpha \in \left( 0,  {\gamma}/{2\eta}\right)$,
then for arbitrary $c$ chosen from the interval
\begin{equation}\label{c_condition}
c\in \left( \frac{4\alpha L q_{\max}}{\eta}\ , \ \frac{4\alpha \mu q_{\min}}{L}-\frac{2\alpha q_{\min}(2L-\mu)}{\eta}
\right),
\end{equation}
there exits a positive constant $0<\delta<1$ such that 
\begin{equation}\label{linear_convg_claim}
\E{\|\bbu^{t+1}-\bbu^{*}\|^2_{\bbG}+c\ \! p^{t+1} \mid \ccalF^t}\leq  (1-\delta)\left(\|\bbu^{t}-\bbu^{*}\|^2_{\bbG}+c\ \! p^t\right).
\end{equation}
\end{theorem}

\begin{proof}
See Appendix \ref{apx_lin_convg_thm}.
\end{proof}

We point out that the linear convergence constant $\delta$ in \eqref{linear_convg_claim} is explicitly available -- see \eqref{linear_convg_70} in Appendix \ref{apx_lin_convg_thm}. It is a function of the strong convexity parameter $\mu$, the Lipschitz continuity constant $L$, lower and upper bounds on the eigenvalues of the matrices $\tbZ$, $\tbZ-\bbZ$, and $\bbI+\bbZ-2\tbZ$, the smallest $q_{\min}$ and largest $q_{\max}$ values for the number of instantaneous functions available at a node, and the stepsize $\alpha$. Insight on the dependence of $\delta$ with problem parameters is offered in Section \ref{sec_convergence_constant}.

The inequality in \eqref{linear_convg_claim} shows that the expected value of the sequence $\|\bbu^{t}-\bbu^{*}\|^2_{\bbG}+c\ \! p^t$ at time $t+1$ given the observation until step $t$ is strictly smaller than the previous iterate at step $t$. Computing the expected value with respect to the initial sigma field $\E{.\mid \ccalF^0}=\E{.}$ implies that in expectation the sequence $\|\bbu^{t}-\bbu^{*}\|^2_{\bbG}+c\ \! p^t$ converges linearly to null, i.e.,
\begin{equation}\label{best_result}
\E{\|\bbu^{t}-\bbu^{*}\|^2_{\bbG}+c\ \! p^{t}}\leq  (1-\delta)^{t} \left(\|\bbu^{0}-\bbu^{*}\|^2_{\bbG}+c\ \! p^0\right).
\end{equation}
We use the result in \eqref{best_result} to establish linear convergence of the sequence of squared norm error $\|\bbx^t-\bbx^*\|^2$ in expectation. 

\begin{corollary}\label{corollary_lin_convg}
Consider the DSA algorithm as defined in \eqref{eqn_sag_update}-\eqref{initial_DSA_local} and recall $\gamma$ is the minimum eigenvalue of the positive definite matrix $\tbZ$. If the hypothesis of Theorem \ref{lin_convg_thm} holds, then there exits a positive constant $0<\delta<1$ such that 
\begin{equation}\label{best222}
\E{\|\bbx^{t}-\bbx^{*}\|^2}\leq  (1-\delta)^{t} \frac{\left(\|\bbu^{0}-\bbu^{*}\|^2_{\bbG}+c\ \! p^0\right).}{\gamma}
\end{equation}
\end{corollary}

\begin{proof}
First note that according to the definitions of $\bbu$ and $\bbG$ in \eqref{z_G_definitions} and the definition of $p^t$ in \eqref{p_definition} , we can write $\|\bbx^{t}-\bbx^{*}\|_{\tbZ}^2\leq \|\bbu^{t}-\bbu^{*}\|^2_{\bbG}+c\ \! p^{t}$. Further, note that the weighted norm $\|\bbx^{t}-\bbx^{*}\|_{\tbZ}^2$ is lower bounded by $\gamma\|\bbx^{t}-\bbx^{*}\|^2$, since $\gamma$ is a lower bound for the eigenvalues of $\tbZ$. Combine these two observations to obtain $\gamma\|\bbx^{t}-\bbx^{*}\|^2\leq \|\bbu^{t}-\bbu^{*}\|^2_{\bbG}+c\ \! p^{t}$. This inequality in conjunction with the expression in \eqref{best_result} follows the claim in \eqref{best222}.
\end{proof}


Corollary \ref{corollary_lin_convg} states that the sequence $\E{\|\bbx^{t}-\bbx^{*}\|^2}$ linearly converges to null. Note that the sequence $\E{\|\bbx^{t}-\bbx^{*}\|^2}$ is not necessarily monotonically decreasing as the sequence $\E{\|\bbu^{t}-\bbu^{*}\|^2_{\bbG}+c\ \! p^t}$ is. The result in \eqref{best222} shows linear convergence of the sequence of variables generated by DSA in expectation. In the following Theorem we show that all local variables $\bbx_n^t$ generated by DSA almost surely converge to the optimal argument of \eqref{opt_problem}.

\begin{theorem}\label{a_s_convg}
Consider the DSA algorithm as defined in \eqref{eqn_sag_update}-\eqref{initial_DSA_local} and assume the same hypothesis of Theorem \ref{lin_convg_thm}.
Then, the sequences of local variables $\bbx_{n}^t$ for all $n=1,\dots, N$ converge almost surely to the optimal argument $\tbx^*$, i.e., 
\begin{equation}\label{a_s_convg_claim}
\lim_{t \to \infty} \bbx_n^t\ = \ \tbx^* \quad a.s.  \qquad \forall\  n=1,\dots, N.
\end{equation}
Further, the almost sure convergence is at least of order $\ccalO(1/t)$.
\end{theorem}

\begin{proof}
See Appendix \ref{apx_a_s_convg}.
\end{proof}


Theorem \ref{a_s_convg} provides almost sure convergence of $\bbx^t$ to the optimal solution $\bbx^*$ which is stronger result than convergence in expectation as in Corollary \ref{corollary_lin_convg}, however, the rate of convergence for the almost sure convergence is  sublinear $\ccalO(1/t)$ which is slower relative to the linear convergence in expectation provided in \eqref{best222}.

\subsection{Convergence constant}\label{sec_convergence_constant}

The constant $\delta$ that controls the speed of convergence can be simplified by selecting specific values for $\eta$, $\alpha$, and $c$. This uncovers connections to the properties of the local objective functions and the network topology. To make this clearer define the condition numbers of the objective function and the graph as 
 \begin{equation}\label{condition_numbers}
 \kappa_f=\frac{L}{\mu}, \qquad \kappa_g=\frac{\max\{\Gamma,\Gamma'\}}{\min\{\gamma,\gamma'\}},
 \end{equation} 
respectively. The condition number of the function is a measure of how difficult it is to minimize the local functions using gradient descent directions. The condition number of the graph is a measure of how slow the graph is in propagating a diffusion process. Both are known to control the speed of convergence of distributed optimization methods. The following corollary illustrates that these condition numbers also determine the convergence speed of DSA.

\begin{corollary}\label{corollary_lin_convg_constant}
Consider the DSA algorithm as defined in \eqref{eqn_sag_update}-\eqref{initial_DSA_local} and assume the same hypothesis of Theorem \ref{lin_convg_thm}. Choose the weight matrices $\bbW$ and $\tbW$ as $\tbW=(\bbI+\bbW)/2$, assign the same number of instantaneous local functions $f_{n,i}$ to each node, i.e., $q_{\min}=q_{\max}=q$, and set the constants $\eta$, $\alpha$ and $c$ as
\begin{equation}\label{constants}
\eta=\frac{2L^2}{\mu}, \quad \alpha=\frac{\gamma\mu}{8L^2},\quad 
c=\frac{q\gamma \mu^2}{4L^3}\left( 1+\frac{\mu}{4L} \right).
\end{equation}
The linear convergence constant $0<\delta<1$ in \eqref{linear_convg_claim} reduces to
\begin{equation}\label{something}
\delta = \min \Bigg[\frac{1}{ 16\kappa_g^2}, \
\frac{1}{q[1+4\kappa_f(1+\gamma/\gamma')]},\
\frac{1}{ 4(\gamma/\gamma')\kappa_f+32\kappa_g\kappa_f^4 }
\Bigg].
\end{equation}
\end{corollary}

\begin{proof}
The given values for $\eta$, $\alpha$, and $c$ satisfy the conditions in Theorem \ref{lin_convg_thm}. Substitute then these values into the expression for $\delta$ in \eqref{linear_convg_70}. Simplify terms and utilize the condition number definitions in \eqref{condition_numbers}. The second term in the minimization in \eqref{linear_convg_70} becomes redundant because it is dominated by the first.
\end{proof}


Observe that while the choices of $\eta$, $\alpha$, and $c$ in \eqref{constants} satisfy all the required conditions of Theorem \ref{lin_convg_thm}, they are not necessarily optimal for maximizing the linear convergence constant $\delta$. Nevertheless, the expression in \eqref{something} shows that the convergence speed of DSA decreases with increases in the graph condition number $\kappa_g$, the local functions condition number $\kappa_f$, and the number of functions assigned to each node $q$. For a cleaner expression observe that both, $\gamma$ and $\gamma'$ are the minimum eigenvalues of the weight matrix $\bbW$ and the weight matrix difference $\tbW-\bbW$. They can therefore be chosen to be of similar order. For reference, say that we choose $\gamma=\gamma'$ so that the ratio $\gamma/\gamma'=1$. In that case, the constant $\delta$ in \eqref{something} reduces to
\begin{equation}\label{approx_constant}
\delta = \min \Bigg[
\frac{1}{16 \kappa_g^2},\
\frac{1}{q(1+8\kappa_f)},\
\frac{1}{4(\kappa_f+8\kappa_f^4\kappa_g)}
\Bigg].
\end{equation}
The three terms in \eqref{approx_constant} establish separate regimes, problems where the graph condition number is large, problems where the number of functions at each node is large, and problems where the condition number of the local functions are large. In the first regime the first term in \eqref{approx_constant} dominates and establishes a dependence in terms of the square of the graph's condition number. In the second regime the middle term dominates and results in an inverse dependence with the number of functions available at each node. In the third regime, the third term dominates. The dependence in this case is inversely proportional to $\kappa_f^4$.

%% file: Simulations.tex

%
\section{Numerical analysis}\label{sec:simulations}
 
We numerically study the performance of the DSA algorithm in solving a logistic regression problem. In this problem we are given $Q=\sum_{n=1}^Nq_n$ training samples that we distribute across $N$ distinct nodes. Denote $q_n$ as the number of samples that are assigned to node $n$. The training points at node $n$ are denoted by $\bbs_{ni}\in \reals^p$ for $i=1,\dots,q_n$ with associated labels $l_{ni}\in \{-1,1\}$. The goal is to predict the probability  $\Pc{l=1\mid \bbs}$ of having label $l=1$ for sample point $\bbs$. The logistic regression model assumes that this probability can be computed as $\Pc{l=1\mid \bbs}=1/(1+\exp(-\bbs^T\bbx))$ given a linear classifier $\bbx$ that is computed based on the training samples. It follows from this model that the regularized maximum log likelihood estimate of the classifier $\bbx$ given the training samples $(\bbs_{ni},l_{ni})$ for $i=1,\ldots,q_n$ and $n=1,\ldots, N$ is the solution of problem
\begin{align}\label{eqn_logistic_regrssion_max_likelihood}
   \tbx^*
          :=   \argmin_{\bbx\in\reals^p}   \frac{\lambda}{2} \|\bbx\|^2 
             +               \sum_{n=1}^N \sum_{i=1}^{q_n} 
                                  \log \Big(1+\exp(-l_{ni}\bbs_{ni}^T\bbx)\Big),
\end{align}
where the regularization term $(\lambda/2)\|\bbx\|^2$ is added to reduce overfitting to the training set. The optimization problem in \eqref{eqn_logistic_regrssion_max_likelihood} can be written in the form of \eqref{opt_problem} by defining the local objective functions $f_{n}$
as 
\begin{equation}\label{logistic_regrssion_local_obj}
   f_n(\bbx) =    \frac{\lambda}{2n} \|\bbx\|^2
                + \sum_{i=1}^{q_n} \log \Big(1+\exp(-l_{ni}\bbs_{ni}^T\bbx)\Big).
\end{equation}
Observe that the local functions $f_{n}$ in \eqref{logistic_regrssion_local_obj} can be written as the average of a set of instantaneous functions $f_{n,i}$ defined as
\begin{equation}\label{logistic_regrssion_local_instan_obj}
   f_{n,i}(\bbx) =    \frac{\lambda}{2n} \|\bbx\|^2
                + {q_n} \log \Big(1+\exp\left(-l_{ni}\bbs_{ni}^T\bbx\right)\Big),  \end{equation}
$\forall\ i= 1,\dots, q_n$. Considering the definitions of instantaneous local functions $f_{n,i}$ in \eqref{logistic_regrssion_local_instan_obj} and local functions $f_{n}$ in \eqref{logistic_regrssion_local_obj}, problem \eqref{eqn_logistic_regrssion_max_likelihood} can be solved using the DSA algorithm.

In our experiments we use a synthetic dataset where components of the feature vectors $\bbs_{{ni}}$ with label $l_{ni}=1$ are generated from a normal distribution with mean $\mu$ and standard deviation $\sigma_+$, while sample points with label $l_{ni}=-1$ are generated from a normal distribution with mean $-\mu$ and standard deviation $\sigma_-$. We consider a network of size $N$ where the edges between nodes are generated randomly with probability $p_c$. The weight matrix $\bbW$ is generated using the Laplacian matrix $\bbL$ of network as
\begin{equation}
\bbW=\bbI-\bbL/\tau,
\end{equation}
where $\tau>(1/2)\lambda_{\max}(\bbL)$. We capture the error of each algorithm by the sum of squared differences of local iterates $\bbx_{n}^t$ from the optimal solution $\tbx^*$ as 
\begin{equation}
e^t=\|\bbx^t-\bbx^*\|^2=\sum_{i=1}^N \|\bbx_i^t-\tbx^*\|^2.
\end{equation}
We use the total number of sample points $Q=500$, feature vectors dimension $p=2$, regularization parameter $\lambda=10^{-4}$, probability of existence of an edge $p_c=0.3$, and $\tau=(2/3)\lambda_{\max}(\bbL)$ . To make the dataset {\it not} linearly separable we set mean to $\mu=2$ and standard deviations to $\sigma_+=\sigma_-=2$. We use a centralized algorithm for computing the optimal argument $\tbx^*$ in all of our experiments.

We provide a comparison of DSA with respect to DGD, EXTRA, stochastic EXTRA, and decentralized SAGA. The stochastic EXTRA is defined by using stochastic gradient in \eqref{sto_gradient} instead of using full gradient as in EXTRA or stochastic averaging gradient as in DSA. The decentralized SAGA is a stochastic version of DGD algorithm that uses stochastic averaging gradient instead of exact gradient which is the naive approach for developing decentralized version of SAGA algorithm.

%
\begin{figure} \centering
\includegraphics[width=0.7\linewidth,height=0.45\linewidth]{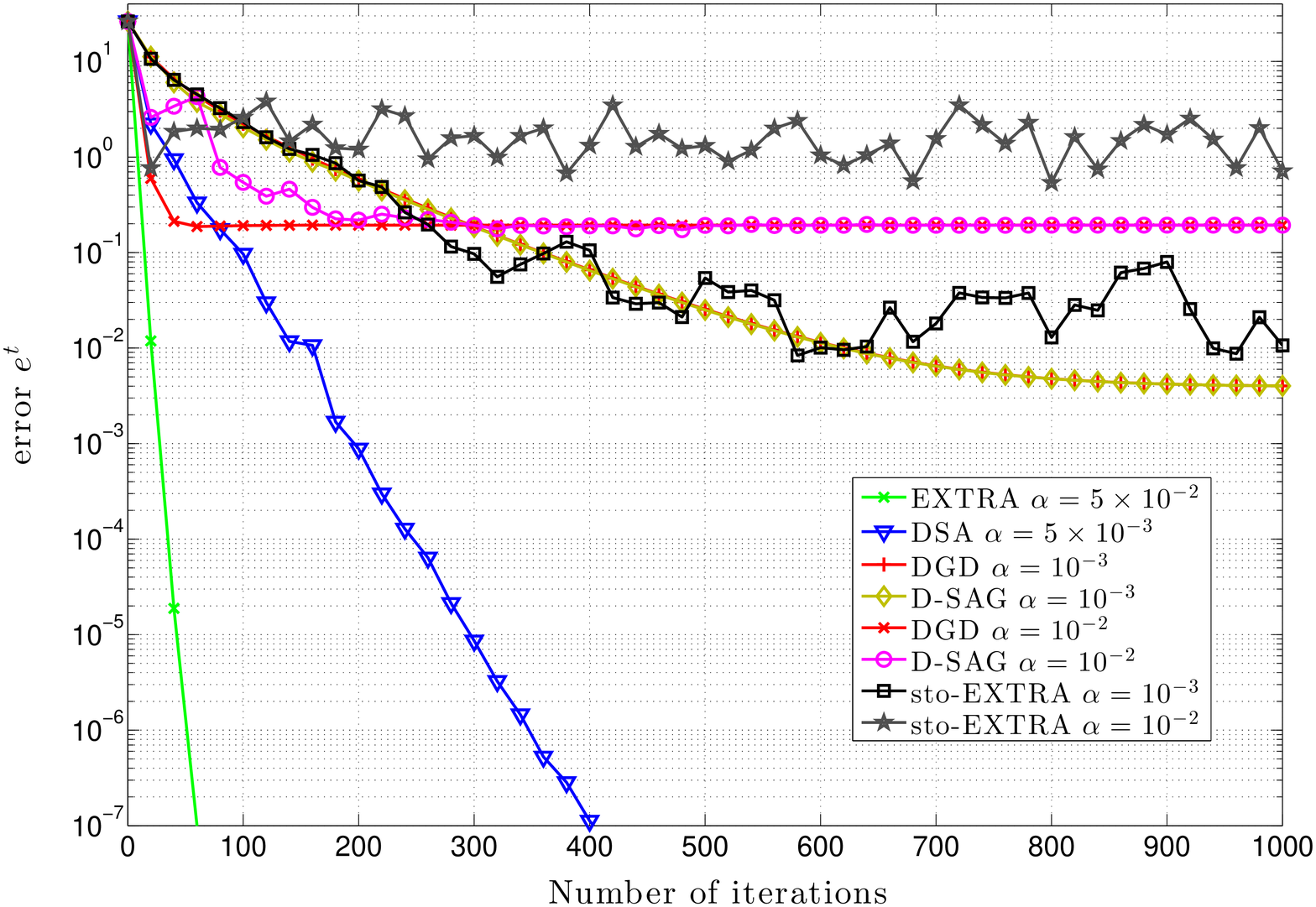}
\caption{ Convergence of DSA, EXTRA, DGD, Stochastic EXTRA, and Decentralized SAGA. Relative distance to optimality $e^t= \|\bbx^t-\bbx^*\|^2$ is shown with respect to number of iterations $t$. DSA and EXTRA converge linearly to the optimal argument $\bbx^*$, while DGD, Stochastic EXTRA, and Decentralized SAGA with constant step sizes converge to a neighborhood of the optimal solution. Smaller choice of stepsize leads to more accurate convergence for these algorithms.
}
\label{fig:ite_comp} \end{figure}

In our experiments the wight matrix $\tbW$ in EXTRA, stochastic EXTRA, and DSA is chosen as $\tbW=(\bbI+\bbW)/2$. Fig. \ref{fig:ite_comp} illustrates the convergence paths of DSA, EXTRA, DGD, Stochastic EXTRA, and Decentralized SAGA with constant step sizes for $N=20$ nodes. For EXTRA and DSA different stepsize are chosen and the best performance for EXTRA and DSA are achieved by $\alpha=5\times10^{-2}$ and $\alpha=5\times10^{-3}$, respectively. As shown in Fig. \ref{fig:ite_comp}, DSA is the only stochastic algorithm that achieves linear convergence. Decentralized SAGA after couple of iterations achieves the performance of DGD and they both can not achieve exact convergence. By choosing smaller stepsize $\alpha=10^{-3}$ they reach more accurate convergence relative to stepsize $\alpha=10^{-2}$, however, the speed of convergence is slower for the smaller stepsize. Stochastic EXTRA also suffers from inexact convergence, but for a different reason. DGD and decentralized SAGA have inexact convergence since they solve a penalty version of the original problem, while stochastic EXTRA can not reach the optimal solution since the noise of stochastic gradient is not vanishing. DSA resolves both issues by combining the idea of stochastic averaging from SAGA to control noise of stochastic gradient and using the double decentralized descent idea of stochastic EXTRA to solve the correct optimization problem. Convergence rate of EXTRA is faster than DSA in terms of number of iterations or equivalently number of communications, however, the complexity of each iteration for EXTRA is higher than DSA. Therefore, it is reasonable to compare performances of these algorithms in terms of number of processed feature vectors. For instance, DSA requires $400$ iterations or equivalently $400$ feature vectors to achieve the error $e^t=10^{-7}$, while to achieve the same accuracy EXTRA requires $60$ iterations which is equivalent to processing $60\times 25=1440$ feature vectors. These numbers show the advantage of DSA relative to EXTRA in requiring less processed feature vectors for achieving a specific accuracy.

%
\begin{figure} \centering
\vspace{-2mm}
\includegraphics[width=0.7\linewidth,height=0.45\linewidth]{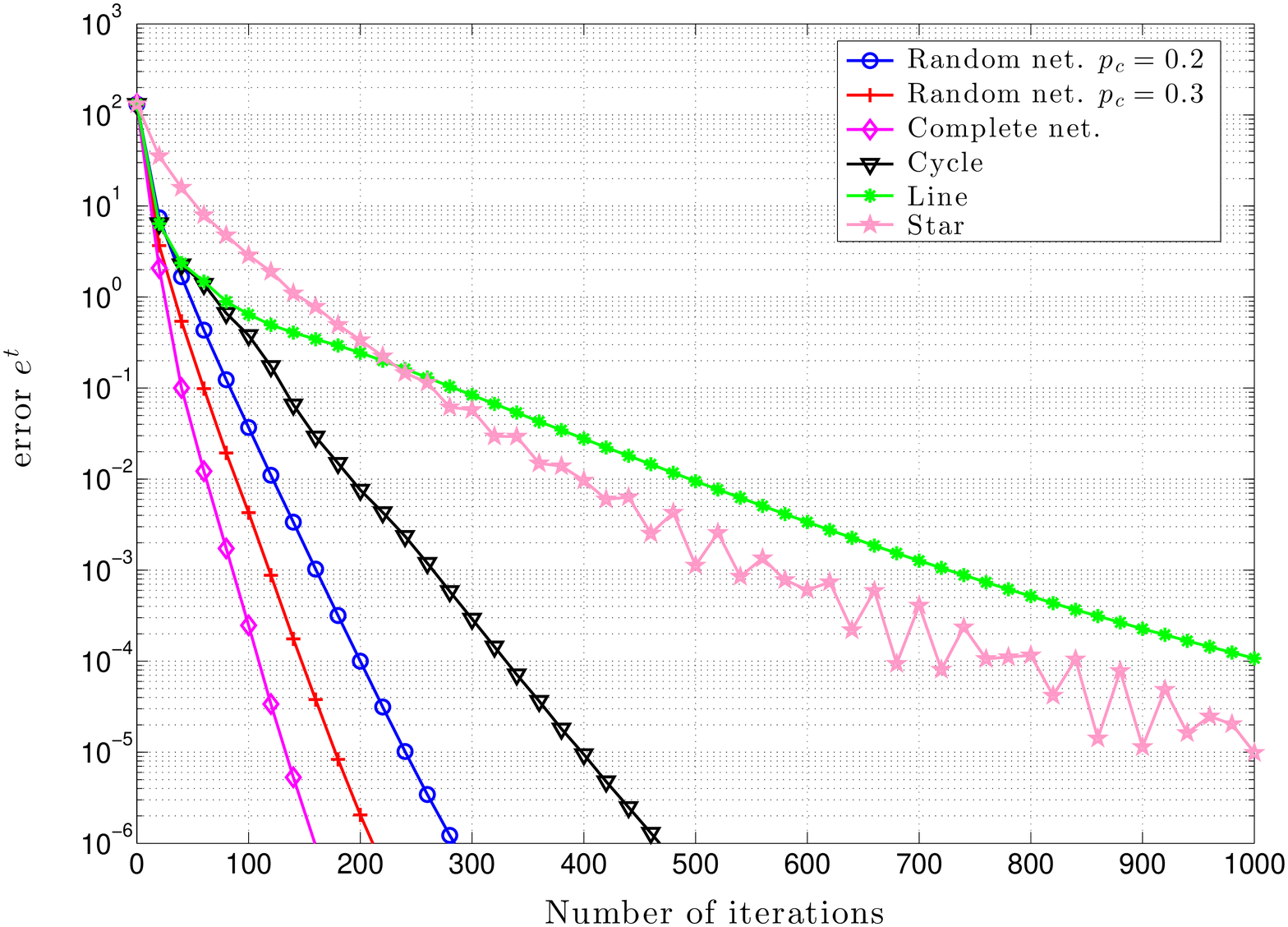}
\caption{ Convergence of DSA for different network topologies. Relative distance to optimality $e^t= \|\bbx_t-\bbx^*\|^2$ is shown with respect to the number iterations $t$. DSA has faster convergence in more connected networks.  
}
\vspace{-2mm}
\label{fig:topology} \end{figure}

We study performances of the DSA algorithm for different topologies. We keep the parameters in Fig. \ref{fig:ite_comp} except we change the size of network to $N=100$ which implies each node has $q_i=5$ sample points. The linear convergence of DSA algorithm for random networks with $p_c=0.2$ and $p_c=0.3$, complete graph, cycle, line and star are shown in Fig. \ref{fig:topology}. As we expect for the topologies that the graph is more connected and the diameter is smaller linear convergence of DSA is faster. The best performance belongs to the complete graph which requires $160$ iterations to achieve the relative error $e^t=10^{-6}$. For random graphs with connectivity probabilities $p_c=0.3$ and $p_c=0.2$ DSA achieves the relative error $e^t=10^{-6}$ after $t=210$ and $t=280$ iterations, respectively. For the cycle graph the number of required iterations for reaching the relative error $e^t=10^{-6}$ is $t=470$, while DSA does not reach this accuracy after $t=1000$ iterations when the graph is a line or star.

%% file: Conclusions.tex

\section{Conclusions}\label{sec_conclusions}

Decentralized double stochastic averaging gradient (DSA) is proposed as an algorithm for solving decentralized optimization problems where the local functions can be written as an average of a set of local instantaneous functions. DSA exploits stochastic averaging gradients in lieu of gradients and mixes information of two consecutive iterates to determine the descent direction. By assuming strongly convex local instantaneous functions with Lipschitz continuous gradients, the DSA algorithm converges linearly to the optimal arguments in expectation. In addition, the sequence of local iterates $\bbx_n^t$ for each node in the network almost surely converges to the optimal argument $\tbx^*$. A comparison between the DSA algorithm and a group of stochastic and deterministic alternatives are provided for solving a logistic regression problem. The numerical results show DSA is the only stochastic decentralized algorithm to reach linear convergence. DSA outperforms decentralized stochastic alternatives in terms of number of required iteration for convergence, and exhibits faster convergence relative to deterministic alternatives in terms of number feature vectors processed until convergence.

%% file: Appendix.tex


\appendix

\section{Proof of Lemma \ref{lemma_gradient_noise_bound}}\label{apx_lemma_grad_noise_bound}


According to the definition of $ \hbg^{t}$ which is the concatenation of local stochastic averaging gradients $ \hbg_n^{t}$ and the fact that expected value of sum is equal to sum of expected values, we can write the expected value $\E{\left\|  \hbg^{t}   - \nabla f(\bbx^*) \right\|^2 \mid\ccalF^t}$ as
\begin{equation}\label{gradient_noise_10}
\E{\left\|  \hbg^{t}    -    \nabla f(\bbx^*) \right\|^2  
\mid\ccalF^t} =
	\sum_{n=1}^N\E{\left\|  \hbg_n^{t} -  \nabla f_n(\tbx^*) \right\|^2  
\mid \ccalF^t}.
\end{equation}
We proceed by finding upper bounds for the summands of \eqref{gradient_noise_10}. 
Observe that using the standard variance decomposition for any random variable vector $\bba$ we can write  $\E{\|\bba\|^2}=\|\E{\bba}\|^2+\E{\|\bba-\E{\bba}\|^2}$. Notice that the same relation holds true when the expectations are computed with respect to a specific field $\ccalF$. By setting $\bba= \hbg_{n}^t  -    \nabla f_n(\tbx^*)  $ and considering the fact that $\E{\bba\mid \ccalF^t}={\nabla} f_{n}(\bbx_n^{t})-\nabla f_n(\tbx^*)$, the variance decomposition implies
\begin{align}\label{gradient_noise_20}
\E{\left\| \hbg_{n}^t  -    \nabla f_n(\tbx^*)   \right\|^2  
\mid \ccalF^t}& =  \left\|  {\nabla} f_{n}(\bbx_n^{t})-\nabla f_n(\tbx^*) \right\|^2
\nonumber \\ 
& \qquad+
\E{\left\| \hbg_{n}^t  -    \nabla f_n(\tbx^*) 
-{\nabla} f_{n}(\bbx_n^{t})+\nabla f_n(\tbx^*)
\right\|^2  
\mid \ccalF^t}.
\end{align}
The next step is to find an upper bound for the last term in \eqref{gradient_noise_20}. Adding and subtracting $ {\nabla} f_{n,i_n^{t}}(\tbx^*)$ and using the inequality $\|\bba+\bbb\|^2\leq 2\|\bba\|^2+2\|\bbb\|^2$ for $\bba= {\nabla} f_{n,i_n^{t}}(\bbx_n^{t}  )- {\nabla} f_{n,i_n^{t}}(\tbx^*)   
-{\nabla} f_{n}(\bbx_n^{t})+\nabla f_n(\tbx^*)$ and $\bbb=-({\nabla} f_{n,i_n^{t}}(\bby_{n,i_n^{t}}^{t}) - {\nabla} f_{n,i_n^{t}}(\tbx^*)   
-({1}/{q_n})\sum_{i=1}^{q_n}{\nabla} f_{n,i}(\bby_{n,i}^{t})   +\nabla f_n(\tbx^*))$ lead to
\begin{align}\label{gradient_noise_30}
&\E{\left\| \hbg_n^t  -    \nabla f_n(\tbx^*) 
-{\nabla} f_{n}(\bbx_n^{t})+\nabla f_n(\tbx^*)
\right\|^2  
\mid \ccalF^t}\\
&\qquad \leq 
2 \E{\left\| {\nabla} f_{n,i_n^{t}}(\bbx_n^{t}  )\!-\! {\nabla} f_{n,i_n^{t}}(\tbx^*)   
\!-\!{\nabla} f_{n}(\bbx_n^{t})\!+\!\nabla f_n(\tbx^*)\right\|^2\! \mid\! \ccalF^t}\nonumber\\
&\qquad \qquad+ 
\!2 \mathbb{E}
\bigg[\Big\| {\nabla} f_{n,i_n^{t}}(\bby_{n,i_n^{t}}^{t})\! -\! {\nabla} f_{n,i_n^{t}}(\tbx^*)
-\!\frac{1}{q_n}\sum_{i=1}^{q_n}{\nabla} f_{n,i}(\bby_{n,i}^{t}) \!  +\!\nabla f_n(\tbx^*)
\Big\|^2\! \mid  \! \ccalF^t\bigg].\nonumber
\end{align}
In this step we use the standard variance decomposition twice to simplify the two expectations in the right hand side of \eqref{gradient_noise_30}. Notice that according to the standard variance decomposition $\E{\|\bba-\E{\bba}\|^2}=\E{\|\bba\|^2}-\|\E{\bba}\|^2$ we obtain $\E{\|\bba-\E{\bba}\|^2}\leq\E{\|\bba\|^2}$. Therefore, by setting $\bby={\nabla} f_{n,i_n^{t}}(\bby_{n,i_n^{t}}^{t}) - {\nabla} f_{n,i_n^{t}}(\tbx^*)   $ and observing that the expected value $\E{{\nabla} f_{n,i_n^{t}}(\bby_{n,i_n^{t}}^{t}) - {\nabla} f_{n,i_n^{t}}(\tbx^*)  \mid \ccalF^{t}}$ is equal to $(1/q_n)\sum_{i=1}^{q_n}{\nabla} f_{n,i}(\bby_{n,i}^{t})   -\nabla f_n(\tbx^*) $ we obtain that 
\begin{align}\label{gradient_noise_50}
&\mathbb{E}\bigg[ \Big\| {\nabla} f_{n,i_n^{t}}(\bby_{n,i_n^{t}}^{t})\! 
-\! {\nabla} f_{n,i_n^{t}}(\tbx^*)   
\!-\!\frac{1}{q_n}\sum_{i=1}^{q_n}{\nabla} f_{n,i}(\bby_{n,i}^{t}) \! +\!\nabla f_n(\tbx^*)
\Big\|^2\mid \ccalF^t\bigg]
\nonumber\\
&\qquad\qquad 
\leq
\E{\left\|{\nabla} f_{n,i_n^{t}}(\bby_{n,i_n^{t}}^{t}) - {\nabla} f_{n,i_n^{t}}(\tbx^*)   \right\|^2\mid \ccalF^t}  .
\end{align}
Moreover, by choosing $\bba={\nabla} f_{n,i_n^{t}}(\bbx_n^{t}  )- {\nabla} f_{n,i_n^{t}}(\tbx^*) $ and noticing the relation for the expected value which is $\E{{\nabla} f_{n,i_n^{t}}(\bbx_n^{t}  )- {\nabla} f_{n,i_n^{t}}(\tbx^*) \mid \ccalF^{t}}={\nabla} f_{n}(\bbx_n^{t})-\nabla f_n(\tbx^*)$, the equality $\E{\|\bba-\E{\bba}\|^2}=\E{\|\bba\|^2}-\|\E{\bba}\|^2$ yields
\begin{align}\label{gradient_noise_60}
\E{\left\|
{\nabla} f_{n,i_n^{t}}(\bbx_n^{t}  )- {\nabla} f_{n,i_n^{t}}(\tbx^*)   
-{\nabla} f_{n}(\bbx_n^{t})+\nabla f_n(\tbx^*)
 \right\|^2\mid \ccalF^t}&=
\E{\left\|{\nabla} f_{n,i_n^{t}}(\bbx_n^{t}  )- {\nabla} f_{n,i_n^{t}}(\tbx^*)    \right\|^2\!\mid \ccalF^t}
\nonumber\\
&\qquad
-\left\|   
{\nabla} f_{n}(\bbx_n^{t})-\nabla f_n(\tbx^*)
\right\|^2.
\end{align}
Substituting the  upper bound in \eqref{gradient_noise_50} and simplification in \eqref{gradient_noise_60} into \eqref{gradient_noise_30}, and considering the expression in  \eqref{gradient_noise_20} lead to
\begin{align}\label{gradient_noise_80}
\E{\left\|  \hbg^{t} - \nabla f(\bbx^*) \right\|^2  \mid \ccalF^t} 
&\leq
	2 \sum_{n=1}^N	\E{\left\|{\nabla} f_{n,i_n^{t}}(\bby_{n,i_n^{t}}^{t}) - {\nabla} f_{n,i_n^{t}}(\tbx^*)   \right\|^2\mid \ccalF^t}- \sum_{n=1}^N  \left\|  {\nabla} f_{n}(\bbx_n^{t})-\nabla f_n(\tbx^*) \right\|^2\nonumber\\
	&\qquad  \quad +
	2 \sum_{n=1}^N \E{\left\|{\nabla} f_{n,i_n^{t}}(\bbx_n^{t}  )
	- {\nabla} f_{n,i_n^{t}}(\tbx^*)    \right\|^2\mid \ccalF^t}.
\end{align}
We proceed by finding an upper bound for the first sum in the right hand side of \eqref{gradient_noise_80}. Notice that if gradients of function $g$ are Lipschitz continuous with parameter $L$, then for any two vectors $\bba_1$ and $\bba_2$ we can write $g(\bba_1)\geq g(\bba_2) + \nabla g(\bba_2)^T(\bba_1-\bba_2)+(1/2L)\|\nabla g(\bba_1) -\nabla g(\bba_2) \|^2$. According to the Lipschitz continuity of instantaneous local functions gradient $\nabla f_{n,i}(\bbx_{n})$, we can write the inequality for $g=f_{n,i}$, $\bba_1=\bby_{n,i}^{t}$ and $\bba_2=\tbx^*$ which is equivalent to 
\begin{align}\label{gradient_noise_90}
\frac{1}{2L}\left\| {\nabla} f_{n,i}(\bby_{n,i}^{t}) - {\nabla} f_{n,i}(\tbx^*) \right\|^2
\leq 
f_{n,i}(\bby_{n,i}^{t}) - f_{n,i}(\tbx^*) 
-\nabla f_{n,i}(\tbx^*)^T (\bby_{n,i}^{t}-\tbx^*). 
\end{align}
Summing up both sides of \eqref{gradient_noise_90} for all $i =1,\dots,q_n$, dividing both sides of the implied inequality by $q_n$ lead to 
\begin{align}\label{gradient_noise_100}
\frac{1}{q_n}\sum_{i=1}^{q_n}\left\| {\nabla} f_{n,i}(\bby_{n,i}^{t}) - {\nabla} f_{n,i}(\tbx^*) \right\|^2
\leq 2L \left[  
\frac{1}{q_n}\sum_{i=1}^{q_n} f_{n,i}(\bby_{n,i}^{t})  -  f_{n,i}(\tbx^*) 
 - \nabla f_{n,i}(\tbx^*)^T (\bby_{n,i}^{t}\! -\! \tbx^*)
  \right].
\end{align}
Since the random functions $f_{n,i_n^t}$ has a uniform distribution over the set $\{f_{n,1},\dots,f_{n,q_n} \}$, we can substitute the left hand side of \eqref{gradient_noise_100} by $	\E{\left\|{\nabla} f_{n,i_n^{t}}(\bby_{n,i_n^{t}}^{t}) - {\nabla} f_{n,i_n^{t}}(\tbx^*)   \right\|^2\mid \ccalF^t}$. Apply this substitution and sum up both sides of \eqref{gradient_noise_100} for $n=1,\dots,N$. According to the definition of sequence $p^t$ in \eqref{p_definition}, if we sum up the right hand side of \eqref{gradient_noise_100} over $n$ it can be simplified as $2Lp^t$. Applying these simplifications we obtain
\begin{equation}\label{gradient_noise_110}
\sum_{n=1}^N\E{\left\|{\nabla} f_{n,\bbtheta_n^{t}}(\bby_n^{t}) - {\nabla} f_{n,\bbtheta_n^{t}}(\tbx^*)   \right\|^2\mid \ccalF^t}
\leq 2Lp^t.
\end{equation}
Substituting the upper bound in \eqref{gradient_noise_110} into \eqref{gradient_noise_80} and simplifying the sum $\sum_{n=1}^N  \left\|  {\nabla} f_{n}(\bbx_n^{t})-\nabla f_n(\tbx^*) \right\|^2$ as $\left\|  {\nabla} f(\bbx^{t})-\nabla f(\bbx^*) \right\|^2$ yield
\begin{align}\label{gradient_noise_120}
\E{\left\|  \hbg^{t}    -   \nabla f(\bbx^*) \right\|^2  
\mid \ccalF^t}
	&\leq
	2\sum_{n=1}^N \E{\left\|{\nabla} f_{n,i_n^{t}}(\bbx_n^{t}  )- {\nabla} f_{n,i_n^{t}}(\tbx^*)    \right\|^2\mid \ccalF^t}	
	-   \left\|  {\nabla} f(\bbx^{t})-\nabla f(\bbx^*) \right\|^2\nonumber\\
	&\qquad
	+4Lp^t.
\end{align}
To show that the sum in the right hand side of \eqref{gradient_noise_120} is bounded above we use the Lipschitz continuity of the instantaneous functions gradients $\nabla f_{n,i}$. Using the same argument from \eqref{gradient_noise_90} to \eqref{gradient_noise_110} we can write 
\begin{align}\label{gradient_noise_130}
&\sum_{n=1}^N \E{\left\|{\nabla} f_{n,i_n^{t}}(\bbx_n^{t}  )- {\nabla} f_{n,i_n^{t}}(\tbx^*)    \right\|^2 \mid  \ccalF^t }
\\
&\qquad \qquad \qquad \qquad
\leq 2L \sum_{n=1}^N
\frac{1}{q_n} \bigg[ 
\sum_{i=1}^{q_n} f_{n,i}(\bbx_{n}^{t}) - f_{n,i}(\tbx^*)  -\nabla f_{n,i}(\tbx^*)^T (\bbx_{n}^{t}-\tbx^*)\bigg].\nonumber
\end{align} 
Considering the definition of the local objective functions $f_{n}(\bbx_n)=(1/q_n) \sum_{i=1}^{q_n} f_{n,i}(\bbx_n)$ and the aggregate function  $f(\bbx):=\sum_{n=1}^N f_{n}(\bbx_n)$, the right hand side of \eqref{gradient_noise_130} can be simplified as
\begin{align}\label{gradient_noise_140}
 \sum_{n=1}^N \E{\left\|{\nabla} f_{n,i_n^{t}}(\bbx_n^{t}  )- {\nabla} f_{n,i_n^{t}}(\tbx^*)    \right\|^2 \mid  \ccalF^t } \leq 2L \left(
 f(\bbx^{t}) - f(\bbx^*)  -\nabla f(\bbx^*)^T (\bbx^{t}-\bbx^*)\right).
\end{align} 
Replacing the sum in \eqref{gradient_noise_120} by the upper bound in \eqref{gradient_noise_140} implies
\begin{align}\label{gradient_noise_150}
\E{\left\|  \hbg^{t}   -  \nabla f(\bbx^*) \right\|^2  
\mid \ccalF^t}
	\leq 
	4Lp^t
	 - \left\|  {\nabla} f(\bbx^{t})-\nabla f(\bbx^*) \right\|^2
	+4L \left(
 f(\bbx^{t}) - f(\bbx^*)  -\nabla f(\bbx^*)^T (\bbx^{t}-\bbx^*)\right). 
\end{align}
Considering the strong convexity of function $f$ with constant $\mu$ we can write
\begin{align}\label{gradient_noise_160}
 \left\|  {\nabla} f(\bbx^{t})-\nabla f(\bbx^*) \right\|^2 \geq 2\mu \left(
 f(\bbx^{t}) - f(\bbx^*)  -\nabla f(\bbx^*)^T (\bbx^{t}-\bbx^*)\right).
\end{align}
Therefore, we can substitute $  \|  {\nabla} f(\bbx^{t})-\nabla f(\bbx^*) \|^2$ in \eqref{gradient_noise_140} by the lower bound in \eqref{gradient_noise_160} and the claim in \eqref{claim_gradient_noise} follows.

\section{Proof of Lemma \ref{lemma_lower_bound_for_z_decrement}}\label{apx_lemma_lower_bound_for_z_decrement}

According to the Lipschitz continuity of the aggregate function gradients $\nabla f(\bbx)$ we can write $(1/L)\|  \nabla f(\bbx^t) -  \nabla f(\bbx^*) \|^2 \leq (\bbx^t-\bbx^*)^T(\nabla f(\bbx^t) -  \nabla f(\bbx^*))$. By adding and subtracting $\bbx^{t+1}$ to the term $\bbx^t-\bbx^*$ and multiplying both sides of the inequality by $2\alpha$ we obtain 
\begin{align}\label{z_decrement_10}
\frac{2\alpha}{L}\!\left\|  \nabla f(\bbx^t)\! -\!  \nabla f(\bbx^*) \right\|^2
	 \leq
	2 \alpha (\bbx^{t+1}-\bbx^*)^T(\nabla f(\bbx^t) -  \nabla f(\bbx^*))
	+2 \alpha (\bbx^t-\bbx^{t+1})^T(\nabla f(\bbx^t)\! -\!  \nabla f(\bbx^*)).
\end{align}
Expanding the difference $\nabla f(\bbx^t) -  \nabla f(\bbx^*)$ as $\hbg^t -  \nabla f(\bbx^*)+\nabla f(\bbx^t)-\hbg^t $ for the first inner product in the right hand side of \eqref{z_decrement_10} implies 
\begin{align}\label{z_decrement_20}
\frac{2\alpha}{L}\!\left\|  \nabla f(\bbx^t)\! -\!  \nabla f(\bbx^*) \right\|^2
	& \leq 
		2 \alpha (\bbx^t-\bbx^{t+1})^T(\nabla f(\bbx^t) -  \nabla f(\bbx^*))
	+2 \alpha (\bbx^{t+1}-\bbx^*)^T(\hbg^t-  \nabla f(\bbx^*))
\nonumber\\
	&\qquad +  2 \alpha (\bbx^{t+1}-\bbx^*)^T(  \nabla f(\bbx^t)-\hbg^{t}).
\end{align}
We proceed to simplify the inner product $2\alpha(\bbx^{t+1}-\bbx^*)^T(\hbg^t-  \nabla f(\bbx^*))$ in the right hand side of \eqref{z_decrement_20} by substituting $\alpha(\hbg^t-  \nabla f(\bbx^*))$ with its equivalent as introduced in \eqref{claim_10}. Applying this substitution the inner product can be simplified as
\begin{align}\label{z_decrement_30}
2 \alpha (\bbx^{t+1}-\bbx^*)^T(\hbg^t-  \nabla f(\bbx^*)) &=
-2\|\bbx^{t+1}-\bbx^*\|_{\bbI+\bbZ-2\tbZ}^2+2 (\bbx^{t+1}-\bbx^*)^T\tbZ(\bbx^{t}-\bbx^{t+1})
\nonumber\\
&\qquad -2 (\bbx^{t+1}-\bbx^*)^T\bbU(\bbv^{t+1}-\bbv^*).
\end{align}
First notice that according to the KKT condition of problem \eqref{constrained_opt_problem} the optimal primal variable satisfies $(\tbZ-\bbZ)^{1/2}\bbx^*=\bb0$ which by considering the definition of matrix $\bbU=(\tbZ-\bbZ)^{1/2}$ we obtain that $\bbU\bbx^{*}=\bb0$. This observation in associations with the update rule of dual variable $\bbv^t$ in \eqref{DSA_dual_update} implies that we can substitute $\bbU(\bbx^{t+1}-\bbx^{*})$ by $\bbv^{t+1}-\bbv^{t}$. Making this substitution into the last summand of the right hand side of \eqref{z_decrement_30} and considering the symmetry of matrix $\bbU$ yield 
\begin{align}\label{z_decrement_40}
2 \alpha (\bbx^{t+1}-\bbx^*)^T(\hbg^t-  \nabla f(\bbx^*))
&=
-2\|\bbx^{t+1}-\bbx^*\|_{\bbI+\bbZ-2\tbZ}^2+2 (\bbx^{t+1}-\bbx^*)^T\tbZ(\bbx^{t}-\bbx^{t+1})
\nonumber\\
&\qquad -2( \bbv^{t+1}-\bbv^{t})^T(\bbv^{t+1}-\bbv^*).
\end{align}
According to the definition of vector $\bbu$ and matrix $\bbG$ in \eqref{z_G_definitions}, the last two summands of \eqref{z_decrement_40} can be simplified as $2(\bbu^{t+1}-\bbu^{t})^T\bbG(\bbu^{*}-\bbu^{t+1})$. 
Moreover, observe that the inner product $2(\bbu^{t+1}-\bbu^{t})^T\bbG(\bbu^{*}-\bbu^{t+1})$ can be simplified as $\|\bbu^t-\bbu^*\|_\bbG^2-\|\bbu^{t+1}-\bbu^*\|_\bbG^2-\|\bbu^{t+1}-\bbu^t\|_\bbG^2$. Applying this simplification into \eqref{z_decrement_40} implies
\begin{align}\label{z_decrement_60}
2 \alpha (\bbx^{t+1}-\bbx^*)^T(\hbg^t-  \nabla f(\bbx^*))
&=
-2\|\bbx^{t+1}-\bbx^*\|_{\bbI+\bbZ-2\tbZ}^2+\|\bbu^t-\bbu^*\|_\bbG^2-\|\bbu^{t+1}-\bbu^*\|_\bbG^2
\nonumber\\
&\qquad
-\|\bbu^{t+1}-\bbu^t\|_\bbG^2.
\end{align}
The next step is to bound above the inner product $2\alpha (\bbx^t-\bbx^{t+1})^T(\nabla f(\bbx^t) -  \nabla f(\bbx^*)) $. Note that for any two vectors $\bba$ and $\bbb$, and any positive scalar $\eta$ the inequality $2\bba^T\bbb\leq \eta\|\bba\|^2+\eta^{-1}\|\bbb\|^2$ holds true. Therefore, by setting $\bba=\bbx^t-\bbx^{t+1}$ and $\bbb=\nabla f(\bbx^t) -  \nabla f(\bbx^*)$ we obtain that 
\begin{align}\label{z_decrement_71}
2\alpha (\bbx^t-\bbx^{t+1})^T(\nabla f(\bbx^t)-  \nabla f(\bbx^*)) 
\leq \frac{\alpha}{\eta} \|\nabla f(\bbx^t) -  \nabla f(\bbx^*)\|^2
+ \alpha\eta  \|\bbx^t-\bbx^{t+1}\|^2.
\end{align}
Now we substitute the terms in the right hand side of \eqref{z_decrement_20} by their simplifications or upper bounds. Replacing the inner product $2 \alpha (\bbx^{t+1}-\bbx^*)^T(\hbg^t-  \nabla f(\bbx^*))$ by the simplification in \eqref{z_decrement_60}, substituting expression $2\alpha (\bbx^t-\bbx^{t+1})^T(\nabla f(\bbx^t) -  \nabla f(\bbx^*))$ by the upper bound in \eqref{z_decrement_71}, and substituting inner product $2 \alpha (\bbx^{t+1}-\bbx^*)^T(  \nabla f(\bbx^t)-\hbg^{t})$ by the sum $2 \alpha (\bbx^{t}-\bbx^*)^T(  \nabla f(\bbx^t)-\hbg^{t})+2 \alpha (\bbx^{t+1}-\bbx^t)^T(  \nabla f(\bbx^t)-\hbg^{t})$ imply
\begin{align}\label{z_decrement_80}
\frac{2\alpha}{L}\left\|  \nabla f(\bbx^t) -  \nabla f(\bbx^*) \right\|^2
	& \leq 
	-2\|\bbx^{t+1}-\bbx^*\|_{\bbI+\bbZ-2\tbZ}^2+\|\bbu^t-\bbu^*\|_\bbG^2-\|\bbu^{t+1}-\bbu^*	\|_\bbG^2 
	\nonumber\\
	&\qquad
	- \|\bbu^{t+1}-\bbu^t\|_\bbG^2+\alpha\eta  \|\bbx^t-\bbx^{t+1}\|^2
	+ \frac{\alpha}{\eta} \|\nabla f(\bbx^t) -  \nabla f(\bbx^*)\|^2
	\nonumber\\
	&\qquad 
	+  2 \alpha (\bbx^{t}-\bbx^*)^T(  \nabla f(\bbx^t)-\hbg^{t})
	+ 2 \alpha (\bbx^{t+1}-\bbx^t)^T(  \nabla f(\bbx^t) - \hbg^{t}).
\end{align}
Considering that $\bbx^t-\bbx^*$ is deterministic given observations until step $t$ and observing the relation $\E{\hbg^{t}\mid \ccalF^t}= \nabla f(\bbx^t) $, we obtain that $\E{ (\bbx^{t}-\bbx^*)^T(  \nabla f(\bbx^t)-\hbg^{t})\mid \ccalF^t}=0$. Therefore, by computing the expected value of both sides of \eqref{z_decrement_80} given the observations until step $t$ and regrouping the terms we obtain 
\begin{align}\label{z_decrement_90}
	\|\bbu^t-\bbu^*\|_\bbG^2-
	\E{\|\bbu^{t+1}-\bbu^*\|_\bbG^2\mid \ccalF^t}
&	\geq
\alpha\left(\frac{2}{L}-\frac{1}{\eta}\right)\left\|  \nabla f(\bbx^t) -  \nabla f(\bbx^*) \right\|^2
\!+\!\E{\|\bbu^{t+1}-\bbu^t\|_\bbG^2\mid\ccalF^t}
\nonumber
\\
&   \quad  
	 +2\E{\|\bbx^{t+1}-\bbx^*\|_{\bbI+\bbZ-2\tbZ}^2\mid \ccalF^t}
	-\alpha\eta \E{ \|\bbx^t-\bbx^{t+1}\|^2\mid\ccalF^t}
\nonumber\\
&\quad - \E{2 \alpha (\bbx^{t+1}-\bbx^t)^T(  \nabla f(\bbx^t)-\hbg^{t})\mid\ccalF^t}.
\end{align}
By applying inequality $2\bba^T\bbb\leq \eta\|\bba\|^2+\eta^{-1}\|\bbb\|^2$ for the choice of vectors $\bba=\bbx^{t+1}-\bbx^t$ and $\bbb= \nabla f(\bbx^t)-\hbg^{t}$, we obtain that  $2 (\bbx^{t+1}-\bbx^t)^T(  \nabla f(\bbx^t)-\hbg^{t})$ is bounded above by $\eta \| \bbx^{t+1}-\bbx^t\|^2+(1/\eta)\|\nabla f(\bbx^t)-\hbg^{t}\|^2$. Replacing $2 (\bbx^{t+1}-\bbx^t)^T(  \nabla f(\bbx^t)-\hbg^{t})$ in \eqref{z_decrement_90} by its upper bound $\eta \| \bbx^{t+1}-\bbx^t\|^2+(1/\eta)\|\nabla f(\bbx^t)-\hbg^{t}\|^2$ yields 
\begin{align}\label{z_decrement_100}
	\|\bbu^t-\bbu^*\|_\bbG^2-
	\E{\|\bbu^{t+1}-\bbu^*\|_\bbG^2\mid \ccalF^t}
&	\geq
\alpha\left(\frac{2}{L}-\frac{1}{\eta}\right)\left\|  \nabla f(\bbx^t) -  \nabla f(\bbx^*) \right\|^2
\!+\!\E{\|\bbu^{t+1}-\bbu^t\|_\bbG^2\mid\ccalF^t}
\nonumber
\\
&   \quad  
	 +2\E{\|\bbx^{t+1}-\bbx^*\|_{\bbI+\bbZ-2\tbZ}^2\mid \ccalF^t}
	\!-\!2\alpha\eta \E{ \|\bbx^t-\bbx^{t+1}\|^2\mid\ccalF^t}
\nonumber\\
&\quad
	-\frac{ \alpha}{\eta}\E{ \| \nabla f(\bbx^t)-\hbg^{t}\|^2\mid\ccalF^t}.
\end{align}
According to the definitions of vector $\bbu$ and matrix $\bbG$ in \eqref{z_G_definitions} the squared norm $\|\bbu^{t+1}-\bbu^t\|_\bbG^2$ can be expanded as $\|\bbx^{t+1}-\bbx^t\|_\tbZ^2+\|\bbv^{t+1}-\bbv^t\|^2$. Making this simplification for $\|\bbu^{t+1}-\bbu^t\|_\bbG^2$ and regrouping the terms in \eqref{z_decrement_100} lead to 
\begin{align}\label{z_decrement_110}
	\|\bbu^t-\bbu^*\|_\bbG^2-
	\E{\|\bbu^{t+1}-\bbu^*\|_\bbG^2\mid \ccalF^t}
&	\geq
\alpha\left(\frac{2}{L}-\frac{1}{\eta}\right)\left\|  \nabla f(\bbx^t) -  \nabla f(\bbx^*) \right\|^2
\\
&\quad 	+\E{\|\bbx^{t+1}-\bbx^t\|_{\tbZ-2\alpha \eta\bbI}^2\mid \ccalF^t}
	+\E{\|\bbv^{t+1}-\bbv^t\|^2\mid \ccalF^t}\nonumber\\
&   \quad  
	 +2\E{\|\bbx^{t+1}\!-\!\bbx^*\|_{\bbI+\bbZ-2\tbZ}^2\!\mid\! \ccalF^t}
	\!-\!\frac{ \alpha}{\eta}\E{ \| \nabla f(\bbx^t)-\hbg^{t}\|^2\!\mid\!\ccalF^t}.\nonumber
\end{align}
We proceed by simplifying $\E{ \| \nabla f(\bbx^t)-\hbg^{t}\|^2\!\mid\!\ccalF^t}$ in \eqref{z_decrement_110}. Note that by adding and subtracting $\nabla f(\bbx^*)$ the expectation can be written as $\E{ \| \nabla f(\bbx^t)-\nabla f(\bbx^*)+\nabla f(\bbx^*)-\hbg^{t}\|^2\mid\ccalF^t}$ and by expanding the squared norm and simplifying the terms we obtain
\begin{align}\label{z_decrement_120}
\E{ \left\| \nabla f(\bbx^t)-\hbg^{t}\right\|^2\mid\ccalF^t}
=\E{ \left\| \hbg^{t}-\nabla f(\bbx^*)\right\|^2\mid\ccalF^t}-\E{ \left\| \nabla f(\bbx^t)- \nabla f(\bbx^*)\right\|^2\mid\ccalF^t}.
\end{align}
Substituting the simplification in \eqref{z_decrement_120} into \eqref{z_decrement_110} yields
\begin{align}\label{z_decrement_130}
	\|\bbu^t-\bbu^*\|_\bbG^2-
	\E{\|\bbu^{t+1}-\bbu^*\|_\bbG^2\mid \ccalF^t}
&	\geq
\frac{2\alpha }{L}\left\|  \nabla f(\bbx^t) -  \nabla f(\bbx^*) \right\|^2
\\
&\quad 	+\E{\|\bbx^{t+1}-\bbx^t\|_{\tbZ-2\alpha\eta\bbI}^2\mid \ccalF^t}
	+\E{\|\bbv^{t+1}-\bbv^t\|^2\mid \ccalF^t}\nonumber\\
&   \quad  
	 +2\E{\|\bbx^{t+1}\!-\!\bbx^*\|_{\bbI+\bbZ-2\tbZ}^2\!\mid\! \ccalF^t}
	\!-\!\frac{\alpha}{\eta}\E{\|\hbg^t-\nabla f(\bbx^*)\|^2 \mid \ccalF^t}.\nonumber
\end{align}
Considering the strong convexity of function $f$ with constant $\mu$ we can write $ \left\|  {\nabla} f(\bbx^{t})-\nabla f(\bbx^*) \right\|^2 \geq 2\mu \left(
 f(\bbx^{t}) - f(\bbx^*)  -\nabla f(\bbx^*)^T (\bbx^{t}-\bbx^*)\right)$. substituting the squared norm $\left\|  {\nabla} f(\bbx^{t})-\nabla f(\bbx^*) \right\|^2$ by this lower bound in \eqref{z_decrement_130} follows
\begin{align}\label{z_decrement_131}
	\|\bbu^t-\bbu^*\|_\bbG^2-
	\E{\|\bbu^{t+1}-\bbu^*\|_\bbG^2\mid \ccalF^t}
&	\geq
\frac{4\alpha\mu }{L}\left(
 f(\bbx^{t}) - f(\bbx^*)  -\nabla f(\bbx^*)^T (\bbx^{t}-\bbx^*)\right)
\\
&\quad 	+\E{\|\bbx^{t+1}-\bbx^t\|_{\tbZ-2\alpha\eta\bbI}^2\mid \ccalF^t}
	+\E{\|\bbv^{t+1}-\bbv^t\|^2\mid \ccalF^t}\nonumber\\
&   \quad  
	 +2\E{\|\bbx^{t+1}\!-\!\bbx^*\|_{\bbI+\bbZ-2\tbZ}^2\!\mid\! \ccalF^t}
	\!-\!\frac{\alpha}{\eta}\E{\|\hbg^t-\nabla f(\bbx^*)\|^2 \mid \ccalF^t}.\nonumber
\end{align}
Substituting the upper bound for the expectation $\E{\| \hbg^t-\nabla f(\bbx^*)\|^2 \mid \ccalF^t}$ in \eqref{claim_gradient_noise} into \eqref{z_decrement_131} and regrouping the terms show validity of the claim in \eqref{claim_z_decrement_prime}.

\section{Proof of Lemma \ref{lemma_p_decrement_2}}\label{apx_lemma_p_decrement_2}

Given the information until time $t$, each auxiliary vector $\bby_{n,i}^{t+1}$ is a random variable that takes values $\bby_{n,i}^{t}$ and $\bbx_{n}^t$ with associated probabilities $1-1/q_n$ and $1/q_n$, respectively. This observation holds since with probability $1/q_n$ node $n$ may choose index $i$ to update at time $t+1$ and with probability $1-(1/q_n)$ choose other indices. Therefore, we can write
\begin{align}\label{p_seq_decrement10}
\E{\frac{1}{q_n}\sum_{i=1}^{q_n}\left(\nabla f_{n,i}(\tbx^*)^T (\bby_{n,i}^{t+1}-\tbx^*)\right)\mid \ccalF^t}  & = \left[1-\frac{1}{q_n}\right]\frac{1}{q_n}\sum_{i=1}^{q_n}\nabla f_{n,i}(\tbx^*)^T \!(\bby_{n,i}^{t}\!-\!\tbx^*)
\nonumber\\
&\qquad
+\frac{1}{q_n} \nabla f_{n}(\tbx^*)^T (\bbx_{n}^{t}-\tbx^*).
\end{align}
Likewise, the distribution of random function $f_{n,i}(\bby_{n,i}^{t+1})$ given observation until time $t$ has two possibilities $f_{n,i}(\bby_{n,i}^{t})$ and $f_{n,i}(\bbx_{n}^{t})$ with associated probabilities $1-1/q_n$ and $1/q_n$, respectively. Hence, we can write $\E{f_{n,i}(\bby_{n,i}^{t+1})\mid \ccalF^t}=(1-1/q_n) f_{n,i}(\bby_{n,i}^{t})+ (1/q_n)f_{n,i}(\bbx_{n}^{t})$. By summing this relation for all $i\in {1,\dots,q_n}$ and divining by $q_n$ we obtain 
\begin{equation}\label{p_seq_decrement20}
\E{\frac{1}{q_n}\!\sum_{i=1}^{q_n} f_{n,i}(\bby_{n,i}^{t+1}) \mid  \ccalF^t}
=
\left[1-\frac{1}{q_n}\right]\frac{1}{q_n}\sum_{i=1}^{q_n} f_{n,i}(\bby_{n,i}^{t}) 
+
\frac{1}{q_n}f_{n}(\bbx_n^t)
\end{equation}
For the simplicity of equations let us define sequence $p_{n}^t$ as
\begin{equation}\label{p_seq_decrement21}
p_n^t:=\frac{1}{q_n}\sum_{i=1}^{q_n} f_{n,i}(\bby_{n,i}^{t})-f_{n}(\tbx^*)-\frac{1}{q_n}\sum_{i=1}^{q_n}\nabla f_{n,i}(\tbx^*)^T (\bby_{n,i}^{t}-\tbx^*).
\end{equation}
Subtracting \eqref{p_seq_decrement10} from \eqref{p_seq_decrement20} and adding $-f_{n}(\tbx^*)$ to the both sides of equality in association with the definition of sequence $p_n^t$ in \eqref{p_seq_decrement21} yield 
\begin{align}\label{p_seq_decrement30}
\E{p_n^{t+1} \mid \ccalF^t}
=
\left[1-\frac{1}{q_n}\right] p_n^t
+\frac{1}{q_n}\left[ f_{n}(\bbx_n^t)-f_{n}(\tbx^*) -\nabla f_{n}(\tbx^*)^T (\bbx_{n}^{t}-\tbx^*)  \right].
\end{align}
We proceed to find and upper bound for the terms in the right hand side of \eqref{p_seq_decrement30}. First note that according to the strong convexity of instantaneous functions $f_{n,i}$ and $f_{n}$ both terms in the right hand side of \eqref{p_seq_decrement30} are non-negative. Observing that the number of instantaneous functions at each node $q_n$ satisfies the condition $q_{\min} \leq q_n \leq q_{\max}$, we obtain 
\begin{equation}\label{p_seq_decrement40}
1-\frac{1}{q_n}\leq 1-\frac{1}{q_{\max}}\ , \qquad \frac{1}{q_n} \leq \frac{1}{q_{\min}}.
\end{equation}
Substituting the upper bounds in \eqref{p_seq_decrement40} into \eqref{p_seq_decrement30}, summing both sides of implied inequality over $n\in\{1,\dots,N\}$, and considering the definitions of optimal argument $\bbx^*=[\tbx^*;\dots;\tbx^*]$ and aggregate function $f(\bbx)=\sum_{n=1}^Nf_{n}(\bbx_n)$ lead to 
\begin{align}\label{p_seq_decrement50}
\sum_{n=1}^N\E{p_n^{t+1} \mid \ccalF^t}
\leq
\left[1-\frac{1}{q_{\max}}\right]\sum_{n=1}^Np_{n}^t+\frac{1}{q_{\min}}\left[ f(\bbx^t)-f(\bbx^*) -\nabla f(\bbx^*)^T(\bbx^{t}-\bbx^*) \right].
\end{align}
Now observe that according to the definitions of sequences $p^t$ and $p_n^t$ in \eqref{p_definition} and \eqref{p_seq_decrement21}, respectively, $p^t$ is the sum of $p_n^t$ for all $n$, i.e. $p^t=\sum_{n=1}^N p_n^t$. Therefore, we can rewrite \eqref{p_seq_decrement50} as
\begin{align}\label{p_seq_decrement60}
\E{
p^{t+1}
 \mid \ccalF^t}
 \leq
 \left[1-\frac{1}{q_{\max}}\right]p^t+\frac{1}{q_{\min}}\left[ f(\bbx^t)-f(\bbx^*) -\nabla f(\bbx^*)^T (\bbx^{t}\!-\!\bbx^*) \right].
\end{align}
Therefore, the claim in \eqref{lemma_p_decrement_claim} is valid.

\section{Proof of Theorem \ref{lin_convg_thm}}\label{apx_lin_convg_thm}

To prove the result of Theorem \ref{lin_convg_thm} first we prove the following Lemma to establish an upper bound for $\|\bbv^t-\bbv^*\|^2.$

\begin{lemma}\label{lemma_upper_bound_for_v_opt_gap}
Consider the DSA algorithm as defined in \eqref{eqn_sag_update}-\eqref{initial_DSA_local}. Further, recall $\gamma'$ as the smallest non-zero eigenvalue and $\Gamma'$ as the largest eigenvalue of matrix $\tbZ-\bbZ$. If Assumptions \ref{ass_wight_matrix_conditions}, \ref{ass_intantaneous_hessian_bounds} and \ref{ass_bounded_stochastic_gradient_norm} hold true, then the squared norm of difference $\|\bbv^t-\bbv^*\|^2$ is bounded above as
\begin{align}\label{claim_v_opt_gap}
\|\bbv^t-\bbv^*\|^2
&\leq
\frac{8}{\gamma'}\E{ \left\|  \bbx^{t+1}\!-\!\bbx^*\right\|_{(\bbI+\bbZ-2\tbZ)^2}^2\!\mid\!\ccalF^t}\!
+\frac{8}{\gamma'}  \E{\left\|\bbx^{t}\!-\!\bbx^{t+1}  \right\|_{\tbZ^2}^2\!\mid\!\ccalF^t}
 +\frac{16\alpha^2L  }{\gamma'} \ p^t
\nonumber\\
&\quad 
	 +\frac{2\Gamma'}{\gamma'}\E{\|\bbv^{t}-\bbv^{t+1}\|^2\mid\ccalF^t}
	+
	 \frac{8\alpha^2 \left(2L- \mu\right)}{\gamma'} \left[ f(\bbx^t)\!-\!f(\bbx^*) \!-\!\nabla f(\bbx^*)^T\! (\bbx^{t}\!-\!\bbx^*) \right].
\end{align}
\end{lemma}

\begin{proof}
Consider the basic inequality $\|\bba+\bbb\|^2\leq 2\|\bba\|^2+2\|\bbb\|^2$ for the case that $\bba=\bbU(\bbv^{t+1}-\bbv^*)$, $\bbb=\bbU(\bbv^{t}-\bbv^{t+1})$  which can be written as
\begin{equation}\label{v_gap_10}
\|\bbU(\bbv^t-\bbv^*)\|^2\leq 2\|\bbU(\bbv^{t+1}-\bbv^*)\|^2+2\|\bbU(\bbv^{t}-\bbv^{t+1})\|^2.
\end{equation}
We proceed by finding an upper bound for $2\|\bbU(\bbv^{t+1}-\bbv^*)\|^2$. Based on the result of Lemma \ref{lemma_important} in \eqref{claim_10}, the term $\bbU(\bbv^{t+1}-\bbv^*)$ is equal to the sum of vectors $\bba+\bbb$ where $\bba= (\bbI+\bbZ-2\tbZ)(\bbx^{t+1}-\bbx^*)-\tbZ(\bbx^{t}-\bbx^{t+1}) $ and $\bbb=-\alpha\hbg^{t} -\nabla f(\bbx^*)$. Therefore, using inequality $\|\bba+\bbb\|^2\leq 2\|\bba\|^2+2\|\bbb\|^2$ we can write 
\begin{align}\label{v_gap_30}
\left\|\bbU(\bbv^{t+1}-\bbv^*)\right\|^2 \leq 2 \left\|  (\bbI+\bbZ-2\tbZ)(\bbx^{t+1}-\bbx^*)-\tbZ(\bbx^{t}\!-\!\bbx^{t+1})  \right\|^2 +2\alpha^2\left\|\hbg^{t} -\nabla f(\bbx^*) \right\|^2.
\end{align}
By using inequality $\|\bba+\bbb\|^2\leq 2\|\bba\|^2+2\|\bbb\|^2$ one more time for vectors $\bba= (\bbI+\bbZ-2\tbZ)(\bbx^{t+1}-\bbx^*)$ and $\bbb=-\tbZ(\bbx^{t}-\bbx^{t+1})  $, we obtain a upper bound for the term $\|  (\bbI+\bbZ-2\tbZ)(\bbx^{t+1}-\bbx^*)-\tbZ(\bbx^{t}-\bbx^{t+1})  \|^2$ and substituting this upper bound into \eqref{v_gap_30} using the definition of weight norm lead to
\begin{align}\label{v_gap_31}
\left\|\bbU(\bbv^{t+1}-\bbv^*)\right\|^2
\leq   4  \left\|  \bbx^{t+1}-\bbx^*\right\|_{(\bbI+\bbZ-2\tbZ)^2}^2
+4\left\|\bbx^{t}-\bbx^{t+1}  \right\|_{\tbZ^2}^2
 +2\alpha^2\left\|\hbg^{t} -\nabla f(\bbx^*) \right\|^2.
\end{align}
Inequality \eqref{v_gap_31} shows an upper bound for $2\|\bbU(\bbv^{t+1}-\bbv^*)\|^2$ in \eqref{v_gap_10}. Moreover, we know that the second term $\|\bbU(\bbv^{t}-\bbv^{t+1})\|^2$ is also bounded above by $\Gamma'\|\bbv^{t}-\bbv^{t+1}\|^2$ where $\Gamma'$ is the largest eigenvalue of matrix $\tbZ-\bbZ=\bbU^2$. Substituting these upper bounds into \eqref{v_gap_10} and computing the expected value of both sides given the information until step $t$ yield 
\begin{align}\label{v_gap_40}
\|\bbU(\bbv^t-\bbv^*)\|^2 
&\leq\ 
8  \E{\left\|  \bbx^{t+1}-\bbx^*\right\|_{(\bbI+\bbZ-2\tbZ)^2}^2\mid\ccalF^t}
+8\E{\left\|\bbx^{t}-\bbx^{t+1}  \right\|_{\tbZ^2}^2\mid\ccalF^t}
\nonumber\\
&\qquad 
 +4\alpha^2\E{\left\|\hbg^{t} -\nabla f(\bbx^*) \right\|^2\mid \ccalF^t}
 +2\Gamma'\E{\|\bbv^{t}-\bbv^{t+1}\|^2\mid \ccalF^t}.
\end{align}
Note that according to the fact that both $\bbv^t$ and $\bbv^*$ lie in the column space of matrix $\bbU$ we obtain $ \|\bbU(\bbv^t-\bbv^*)\|^2\geq\gamma' \|\bbv^t-\bbv^*\|^2$. Substituting this lower bound for $ \|\bbU(\bbv^t-\bbv^*)\|^2$ in \eqref{v_gap_40} and multiplying both sides of the imposed inequality by $\gamma'$ yield 
\begin{align}\label{v_gap_50}
\|\bbv^t-\bbv^*\|^2
&\leq
\frac{8}{\gamma'} \E{ \left\|  \bbx^{t+1}-\bbx^*\right\|_{(\bbI+\bbZ-2\tbZ)^2}^2\mid\ccalF^t}
+  \frac{8}{\gamma'}\E{\left\|\bbx^{t}-\bbx^{t+1}  \right\|_{\tbZ^2}^2\mid\ccalF^t}
\nonumber\\
&\qquad 
 +\frac{4\alpha^2}{\gamma'} \E{\left\|\hbg^{t} -\nabla f(\bbx^*) \right\|^2\mid\ccalF^t}
 +\frac{2\Gamma'}{\gamma'}\E{\|\bbv^{t}-\bbv^{t+1}\|^2\mid\ccalF^t}.
\end{align}
Substituting $\E{\|\hbg^{t} -\nabla f(\bbx^*) \|^2\mid\ccalF^t}$ in the right hand side of \eqref{v_gap_50} by its upper bound in \eqref{claim_gradient_noise} follows the claim in \eqref{claim_v_opt_gap}.
\end{proof}


Using the result in Lemma \ref{lemma_upper_bound_for_v_opt_gap} we show linear convergence of the sequence $\|\bbu^t-\bbu^*\|_\bbG^2+c\ p^t$ as follows.

\textbf{Proof of Theorem \ref{lin_convg_thm}:}
Proving the linear convergence claim in \eqref{linear_convg_claim} is equivalent to showing that 
\begin{align}\label{linear_convg_20}
\delta\|\bbu^{t}-\bbu^*\|_\bbG^2+\delta c\ p^{t}
\leq 
\|\bbu^{t}-\bbu^*\|_\bbG^2-\E{\|\bbu^{t+1}-\bbu^*\|_\bbG^2\mid \ccalF^t}
+c\ (p^t-\E{p^{t+1}\mid \ccalF^t}).
\end{align}
Substituting the terms $\E{\|\bbu^{t+1}-\bbu^*\|_\bbG^2\mid \ccalF^t}$ and $\E{p^{t+1}\mid \ccalF^t}$ by their upper bounds as introduced in Lemma \ref{lemma_lower_bound_for_z_decrement} and Lemma \ref{lemma_p_decrement_2}, respectively, yield a sufficient condition for the claim in \eqref{linear_convg_20} as 
\begin{align}\label{linear_convg_30}
\delta\|\bbu^{t}-\bbu^*\|_\bbG^2+\delta c\ p^{t}
&\leq 
\E{\|\bbx^{t+1}-\bbx^t\|_{\tbZ-2\alpha\eta\bbI}^2\mid \ccalF^t}
	+\E{\|\bbv^{t+1}-\bbv^t\|^2\mid \ccalF^t}\nonumber\\
&   \quad  
	 +2\E{\|\bbx^{t+1}-\bbx^*\|_{\bbI+\bbZ-2\tbZ}^2\mid \ccalF^t}
	+\left(\frac{c}{q_{\max}}-\frac{4\alpha L}{\eta}\right)p^t\nonumber\\
	&\quad 
	+\left[\frac{4\alpha\mu }{L} -   \frac{2\alpha( 2L-\mu)}{\eta}-\frac{c}{q_{\min}}  \right]
	\left[ f(\bbx^t)-f(\bbx^*) -\nabla f(\bbx^*)^T (\bbx^{t}-\bbx^*) \right].
\end{align}
We emphasize that if inequality \eqref{linear_convg_30} holds then the inequalities in \eqref{linear_convg_20} and \eqref{linear_convg_claim} are valid. Note that $\|\bbu^{t}-\bbu^*\|_\bbG^2$ in the left hand side of \eqref{linear_convg_30} can be simplified as $\|\bbx^t-\bbx^*\|_\tbZ^2+\|\bbv^t-\bbv^*\|^2$. Considering the definition of $\Gamma$ as the maximum eigenvalue of matrix $\tbZ$, we can conclude that $\|\bbx^t-\bbx^*\|_\tbZ^2$ is bounded above by $\Gamma\|\bbx^t-\bbx^*\|^2$. Considering this relation and observing the upper bound for $\|\bbv^t-\bbv^*\|^2$ in \eqref{claim_v_opt_gap}, we obtain that $\|\bbu^{t}-\bbu^*\|_\bbG^2= \|\bbx^t-\bbx^*\|_\tbZ^2+\|\bbv^t-\bbv^*\|^2$ is bounded above as
\begin{align}\label{linear_convg_40}
\|\bbu^{t}-\bbu^*\|_\bbG^2
&\ \leq\
\frac{8}{\gamma'}\E{ \left\|  \bbx^{t+1}-\bbx^*\right\|_{(\bbI+\bbZ-2\tbZ)^2}^2\mid\ccalF^t}
+\frac{8}{\gamma'}  \E{\left\|\bbx^{t}\!-\!\bbx^{t+1}  \right\|_{\tbZ^2}^2\!\mid\!\ccalF^t}
 +\frac{16\alpha^2L  }{\gamma'} \ p^t
\nonumber\\
&\qquad 
	 +\frac{2\Gamma'}{\gamma'}\E{\|\bbv^{t}-\bbv^{t+1}\|^2\mid\ccalF^t}+\Gamma \|\bbx^t-\bbx^*\|^2
	 \nonumber\\
&\qquad 
	+
	 \frac{8\alpha^2 \left(2L- \mu\right)}{\gamma'} \left[ f(\bbx^t)-f(\bbx^*) -\nabla f(\bbx^*)^T (\bbx^{t}-\bbx^*) \right].
\end{align}
Further, substitute the squared norm $\|\bbx^{t}-\bbx^*\|^2$ by the upper bound $(2/\mu)(f(\bbx^t)-f(\bbx^*) -\nabla f(\bbx^*)^T (\bbx^{t}-\bbx^*) ) $ to obtain
\begin{align}\label{linear_convg_41}
\|\bbu^{t}-\bbu^*\|_\bbG^2
&
\leq
\frac{8}{\gamma'}\E{ \left\|  \bbx^{t+1}-\bbx^*\right\|_{(\bbI+\bbZ-2\tbZ)^2}^2\mid\ccalF^t}
+\frac{8}{\gamma'}  \E{\left\|\bbx^{t}-\bbx^{t+1}  \right\|_{\tbZ^2}^2\mid\ccalF^t}
 +\frac{16\alpha^2L  }{\gamma'} \ p^t
\nonumber\\
&\qquad 
	 +\frac{2\Gamma'}{\gamma'}\E{\|\bbv^{t}-\bbv^{t+1}\|^2\mid\ccalF^t}
	 \nonumber\\
&\qquad
	+
	\left( \frac{8\alpha^2 \left(2L- \mu\right)}{\gamma'} +\frac{2\Gamma}{\mu} \right) \left[ f(\bbx^t)-f(\bbx^*) -\nabla f(\bbx^*)^T(\bbx^{t}-\bbx^*) \right].
\end{align}
Replacing $\|\bbu^{t}-\bbu^*\|_\bbG^2$ in \eqref{linear_convg_30} by the upper bound \eqref{linear_convg_41} and regrouping the terms lead to 
\begin{align}\label{linear_convg_50}
0
&\leq 
\E{\| \bbx^{t+1}-\bbx^{t}\|_{\tbZ-\alpha(\eta+\eta)\bbI-\frac{8\delta}{\gamma'} \tbZ^2}^2 \mid \ccalF^t}
+\E{\| \bbx^{t+1}-\bbx^{*}\|_{(\bbI+\bbZ-2\tbZ)^{\frac{1}{2}} \left[2\bbI-\frac{8\delta } {\gamma'}(\bbI+\bbZ-2\tbZ) \right] (\bbI+\bbZ-2\tbZ)^{\frac{1}{2}}  }^2\! \mid \!\ccalF^t}
\nonumber\\
&\qquad +\E{\| \bbv^{t+1}\!\!-\!\bbv^{t}\|_{(1\!-\!\frac{2\delta\Gamma'}{\gamma'})\bbI}^2\! \mid\! \ccalF^t}
\!+\!\left[\frac{c}{q_{\max}}\!-\!\frac{4\alpha L}{\eta}
\!-\!\delta c\! -\!\frac{16\delta\alpha^2L  }{\gamma'}  \right]p^t
\nonumber\\
&\qquad 
+\bigg[\frac{4\alpha\mu}{L}-\frac{2\alpha (2L-\mu)}{\eta} -\!\frac{c}{q_{\min}} - \frac{8\delta\alpha^2 \left(2L- \mu\right)}{\gamma'} -\frac{2\delta\Gamma}{\mu}\bigg]
(f(\bbx^t)\!-\!f(\bbx^*) \!-\!\nabla f(\bbx^*)^T\! (\bbx^{t}\!-\!\bbx^*) ) .
\end{align}
Notice that if the inequality in \eqref{linear_convg_50} holds true, then the relation in \eqref{linear_convg_30} is valid and as we mentioned before the claim in \eqref{linear_convg_20} holds. 
To verify the sum in the right hand side of \eqref{linear_convg_50} is always positive and the inequality is valid, we enforce each summands in the right hand side of \eqref{linear_convg_50} to be non-negative. Therefore, the following conditions should be satisfied
\begin{align}\label{linear_convg_60}
&\gamma-\alpha(\eta+\eta)- \frac{8\delta}{\gamma'}\Gamma^2\geq0 , \quad 
2-\frac{8\delta } {\gamma'} \lambda_{max}(\bbI+\bbZ-2\tbZ) \geq0, \quad  1-\frac{2\delta\Gamma'}{\gamma'}\geq0, 
  \nonumber\\
&
\frac{c}{q_{\max}}-\frac{4\alpha L}{\eta}
-\delta c -\frac{16\delta\alpha^2L }{\gamma'}\geq0, \quad   
\frac{4\alpha\mu}{L}-\frac{2\alpha (2L-\mu)}{\eta} -\!\frac{c}{q_{\min}}
- \frac{8\delta\alpha^2 \left(2L- \mu\right)}{\gamma'} -\frac{2\delta\Gamma}{\mu} 
\geq0.
\end{align}
Recall that $\gamma$ is the smallest eigenvalue of positive definite matrix $\bbZ$.
All the inequalities in \eqref{linear_convg_60} are satisfied, if $\delta $ is chosen as 
\begin{align}\label{linear_convg_70}
&\delta=\min \Bigg\{\frac{(\gamma-2\alpha\eta)\gamma'}{ 8\Gamma^2},
\frac{\gamma'}{4\lambda_{max}(\bbI+\bbZ-2\tbZ)},
\frac{\gamma'}{2\Gamma'}, \frac{\gamma'(c\eta-4\alpha Lq_{\max})}{\eta q_{\max}(c\gamma'+16\alpha^2L)},
\nonumber\\
& 
\qquad \qquad \qquad \qquad \qquad \qquad 
\left[\frac{4\alpha\mu}{L}-\frac{2\alpha (2L-\mu)}{\eta} -\!\frac{c}{q_{\min}}  \right]
\left[\frac{8\alpha^2 \left(2L- \mu\right)}{\gamma'} 
\!+\frac{2\Gamma}{\mu} \right]^{-1}
\Bigg\}.
\end{align}
where $\eta$, $c$ and $\alpha$ are selected from the intervals 
\begin{align}\label{linear_convg_80}
\eta\in\left(
 \frac{L^2q_{\max}}{\mu q_{\min} }+\frac{L^2}{\mu}-\frac{L}{2}\ , \ \infty
\right),\ \alpha \in \left( 0\ , \ \frac{\gamma}{2\eta}\right),\ c\in \left( \frac{4\alpha L q_{\max}}{\eta}\ , \ \frac{4\alpha \mu q_{\min}}{L}-\frac{2\alpha q_{\min}(2L-\mu)}{\eta}
\right).
\end{align}
Notice that considering the conditions for the variables $\eta$, $\alpha$ and $c$ in \eqref{linear_convg_80}, the constant $\delta$ in \eqref{linear_convg_70} is strictly positive $\delta>0$. Moreover, according to the definition in \eqref{linear_convg_70} the constant $\delta$ is smaller than $\gamma'/2\Gamma'$ which leads to the conclusion that $\delta\leq1/2<1$. Therefore, we obtain that $0<\delta<1$ and the claim in \eqref{linear_convg_claim} is valid.

\section{Proof of Theorem \ref{a_s_convg}}\label{apx_a_s_convg}

The proof uses the relationship in the statement \eqref{linear_convg_claim} of Theorem \ref{lin_convg_thm} to build a supermartingale sequence. To do this define the stochastic processes $\zeta^t$ and $\beta^t$ as
\begin{equation}\label{a_s_convg_10}
\zeta^t:= \|\bbu^{t}-\bbu^{*}\|^2_{\bbG}+c\ \! p^t, \quad \beta^t:= \delta\left(\|\bbu^{t}-\bbu^{*}\|^2_{\bbG}+c\ \! p^t\right).
\end{equation}
Note that the stochastic processes $\zeta^t$ and $\beta^t$ are alway non-negative. Let now $\ccalF_t$ be a sigma-algebra measuring $\zeta^t$, $\beta^t$, and $\bbu^t$. Considering the definitions of $\zeta^t$ and $\beta^t$ and the relation in  \eqref{linear_convg_claim} we can write 
\begin{align}\label{a_s_convg_20}
\E{\zeta^{t+1} \mid \ccalF^t} \leq \zeta^{t}-  \beta^t.
\end{align}
Since the sequences $\alpha^t$ and $\beta^t$ are nonnegative it follows from \eqref{a_s_convg_20} that they satisfy the conditions of the supermartingale convergence theorem -- {see e.g. theorem E$7.4$ \cite{Solo}} . Therefore, we obtain that: (i) The sequence $\zeta^t$ converges almost surely. (ii) The sum $\sum_{t=0}^{\infty}\beta^t < \infty$ is almost surely finite. The definition of $\beta^t$ in \eqref{a_s_convg_10} implies that 
\begin{equation}\label{a_s_convg_30}
   \sum_{t=0}^{\infty} \delta\left(\|\bbu^{t}-\bbu^{*}\|^2_{\bbG}+c\ \! p^t\right)< \infty, 
       \qquad\text{a.s.}
\end{equation}
Since $\|\bbx^t-\bbx^*\|^2_{\tbZ} \leq \|\bbu^{t}-\bbu^{*}\|^2_{\bbG}+c\ \! p^t$ and the eigenvalues of $\tbZ$ are lower bounded by $\gamma$ we can write $\gamma\|\bbx^t-\bbx^*\|^2 \leq \|\bbu^{t}-\bbu^{*}\|^2_{\bbG}+c\ \! p^t$. This inequality in association with the fact that the sum in \eqref{a_s_convg_30} is finite leads to 
\begin{equation}\label{a_s_convg_40}
   \sum_{t=0}^{\infty} \delta\ \!\gamma\ \! \|\bbx^t-\bbx^*\|^2< \infty, 
       \qquad\text{a.s.}
\end{equation}
Observing the fact that $\delta $ and $\gamma$ are positive constants, we can conclude from \eqref{a_s_convg_40} that the sequence $\|\bbx^t-\bbx^*\|^2$ is almost surely summable and the it converges with probability 1 to null at least in the order of $\ccalO(1/t)$. Almost sure convergence of sequence to null follows the claim in \eqref{a_s_convg_claim}.